\documentclass{elsarticle}

\usepackage{moreverb}
\usepackage[a4paper, total={6in, 8in}]{geometry}

\makeatletter
\def\ps@pprintTitle{%
  \let\@oddhead\@empty
  \let\@evenhead\@empty
 \def\@oddfoot{}%
 \let\@evenfoot\@oddfoot}
\makeatother

\usepackage[utf8x]{inputenc}

\usepackage{microtype}
\usepackage{graphicx}
\usepackage{subcaption}
\usepackage{booktabs} 
\usepackage{amsmath}
\usepackage{enumerate}

\usepackage{tikz}

\usepackage{algorithm}
\usepackage{algorithmicx}
\usepackage[noend]{algpseudocode}

\usepackage{soul}
\usepackage{amsfonts}

\usepackage{amsthm}
\usepackage{graphicx}
\usepackage{bm,bbm}
\usepackage{comment}
\usepackage{rotating}
\usepackage{color}
\usepackage{overpic}
\usepackage{multirow}
\usepackage{multicol}

\usepackage{soul}
\usepackage{numprint}
\usepackage{amsmath}

\newtheorem{theorem}{Theorem}[section]
\newtheorem{lemma}[theorem]{Lemma}

\newtheorem{proposition}[theorem]{Proposition}

\newtheorem{definition}[theorem]{Definition}

\newtheorem{example}[theorem]{Example}

\newcommand{\TEXTFONTB}[1]{#1}
\newcommand{\va}{\TEXTFONTB{a}}

\newcommand{\ve}{\TEXTFONTB{e}}
\newcommand{\vf}{\TEXTFONTB{f}}

\newcommand{\vh}{\TEXTFONTB{h}}

\newcommand{\vw}{\TEXTFONTB{w}}

\newcommand{\TEXTFONTM}[1]{#1}
\newcommand{\mA}{\TEXTFONTM{A}}
\newcommand{\mB}{\TEXTFONTM{B}}
\newcommand{\mC}{\TEXTFONTM{C}}

\newcommand{\mF}{\TEXTFONTM{F}}

\newcommand{\mH}{\TEXTFONTM{H}}
\newcommand{\mI}{\TEXTFONTM{I}}

\newcommand{\mP}{\TEXTFONTM{P}}

\newcommand{\mV}{\TEXTFONTM{V}}
\newcommand{\mW}{\TEXTFONTM{W}}

\newcommand{\sN}{N}

\newcommand{\sU}{U}

\newcommand{\TEXTFONTT}[1]{\mathcal{#1}}

\newcommand{\tD}{\TEXTFONTT{D}}
\newcommand{\tE}{\TEXTFONTT{E}}

\newcommand{\tG}{\TEXTFONTT{G}}
\newcommand{\tH}{\TEXTFONTT{H}}
\newcommand{\tI}{\TEXTFONTT{I}}

\newcommand{\tX}{\TEXTFONTT{X}}
\newcommand{\tY}{\TEXTFONTT{Y}}

\newcommand{\cC}{{\cal C}}

\newcommand{\normF} [1]{\big\lVert#1\big\rVert_{\mathcal{F}}}

\newcommand{\REAL}{\mathbbm{R}}

\newcommand{\ccirc}{\otimes}

\newcommand{\rank} {\textrm{rank}}

\definecolor{MyDarkGreen}{rgb}{0,0.45,0}

\def\trait #1 #2 #3 {\vrule width #1pt height #2pt depth #3pt}
\def\fin{\hfill
        \trait .3 5 0
        \trait 5 .3 0
        \kern-5pt
        \trait 5 5 -4.7
        \trait 0.3 5 0
\medskip}
\newcommand{\ENDPROOF}{\fin}

\newcommand{\col}{\textsf{col}}

\begin{document}

\begin{frontmatter}

  \title{Nonnegative Canonical Tensor Decomposition with Linear
    Constraints: nnCANDELINC}
  
  \author[1]{Boian Alexandrov}
  \author[2]{Derek DeSantis}
  \author[3]{Gianmarco Manzini}
  \author[4]{Erik Skau}

  \address[1]{Group~T-1, Theoretical Division, Los Alamos
    National~Laboratory, New Mexico, USA}

  \address[2]{Group~CCS-2, Computer, Computational, and Statistical
    Sciences Division,
    Los Alamos\\ National~Laboratory, New Mexico,
    USA}

  \address[3]{Group~T-5, Theoretical Division, Los Alamos
    National~Laboratory, New Mexico, USA}

  \address[4]{Group~CCS-3, Computer, Computational, and Statistical
    Sciences Division,
    Los Alamos\\ National~Laboratory, New Mexico,
    USA}

  \begin{abstract}
    There is an emerging interest in tensor factorization applications
    in big-data analytics and machine learning. To speed up the
    factorization of extra-large datasets, organized in
    multidimensional arrays (aka tensors), easy to compute
    compression-based tensor representations, such as Tucker and
    Tensor Train formats, are used to approximate the initial
    large-tensor. Further, tensor factorization is used to extract
    latent features that can facilitate discoveries of new mechanisms
    and signatures hidden in the data, where the explainability of the
    latent features is of principal importance. Nonnegative tensor
    factorization extracts latent features that are naturally sparse
    and parts of the data, which makes them easily
    interpretable. However, to take into account available domain
    knowledge and subject matter expertise, additional constraints
    often need to be imposed, which lead us to Canonical decomposition
    with linear constraints (CANDELINC), a Canonical Polyadic
    Decomposition with rank deficient factors. In CANDELINC, Tucker
    compression is used as a pre-processing step, which leads to a
    larger residual error but to more explainable latent
    features. Here, we propose a nonnegative CANDELINC (nnCANDELINC)
    accomplished via a specific nonnegative Tucker decomposition; we
    refer to as minimal or canonical nonnegative Tucker. We derive
    several results required to understand the specificity of
    nnCANDELINC, focusing on the difficulties of preserving the
    nonnegative rank to its Tucker core and comparing the real-valued
    to the nonnegative case. Finally, we demonstrate nnCANDELINC
    performance on synthetic and real-world examples.
  \end{abstract}
  
  
\end{frontmatter}





\newcommand{\EOD}{\end{document}}
\newcommand{\NOTE}[1]{~{\textbf{\big[#1\big]}}}

\newcommand{\matricization}{\operatorname{unfold}}
\newcommand{\multirank}{\operatorname{\mu \rank}}
\newcommand{\cone}{\operatorname{cone}}
\newcommand{\rmkes}[1]{%
{\textcolor{red}{[Erik: \small #1]}}%
}
\newcommand{\rmkdd}[1]{%
{\textcolor{blue}{[Derek: \small #1]}}%
}
\setcounter{MaxMatrixCols}{20}

\section{Introduction}

Large amounts of high-dimensional data are constantly generated by
sensor networks; large-scale scientific experiments; massive computer
simulations; complex engineering activities; electronic
communications; social networks, and many other sources
\cite{Big-Data}.
Utilizing such big-data for decision making, emergency response, and
data-driven science requires understanding the processes underlying
the data \cite{franke2016statistical}.
High-dimensional data are naturally organized in tensors (i.e.,
multi-dimensional arrays).
Tensor factorization is a cutting-edge factor analysis that can serve
for latent features extraction, dimensional reduction, blind source
separation, data mining, pattern recognition, subspace learning, data
fusion, compression, and many other
applications~\cite{Kolda-Bader:2009,cichocki2017tensor}.
A tensor factorization's main objective is to decompose
high-dimensional data into factor matrices and one, or in the case of
tensor networks \cite{oseledets2011tensor}, several core-tensors of a
smaller size.

The number of the tensor entries scales exponentially with tensor
dimension, which leads to exponential scaling of the burden of any
tensor computation, in terms of storage and floating point
operations. This phenomenon is known as the \emph{curse of
dimensionality}. One way to speed up tensor calculations and decrease
the needed storage is to use stable compression-based representations
of the large initial tensors, and then to extract the needed
information from the compressed data
\cite{vervliet2018compressed}. Some of the proposed stable
compression-based formats are Tucker \cite{Tucker:1966} (related to
the multirank of a tensor \cite{Hackbusch:2012}) and Tensor Train (TT)
formats \cite{oseledets2011tensor}, which need $O(dnr+r^d)$ and $O(dnr
+ (d − 2)r^3)$ parameters, respectively, vs. $O(n^d)$ entries of the
full tensor (here $d$ is the tensor dimension, $n$ is the number of
entries in each dimension, and $r$ is the Tucker/TT ranks used in
compression). Canonical Polyadic Decomposition
~\cite{Hitchcock:1927,Harshman:1970}(CPD), related to the rank of the
tensor \cite{Hackbusch:2012}, also offers a good compression, however,
computation of the tensor rank is an NP-hard problem
\cite{haastad1989tensor}, and ill-conditioned decompositions and
ill-posed optimization problems often remain unsolved
\cite{de2008tensor}.

Another problem is that the existing datasets are formed by directly
observable quantities, while the underlying processes (features or
variables) usually remain unobserved, hidden, or latent
\cite{everett2013introduction}.
This necessitates the ability to identify and extract
\emph{explainable} latent features needed to identify essential
signatures that are manifestation of the processes and causalities
hidden in large high-dimensional datasets. Imposing various
constraints on the factors, reflecting available prior information,
usually helps to mitigate this problem.

Many types of real-world data (e.g., density, energy, spectral power,
population, pixels, probabilities, frequencies of appearance, etc.)
are naturally nonnegative and the extracted features will lose their
meaning if the nonnegativity is not preserved.  Tensor factorizations
with nonnegative constraint extract nonnegative latent features formed
by only positive combinations, which favors parts based sparse
representation where extracted features are parts of the original data
\cite{lee1999learning}.
Importantly, because the extracted features are parts of the original
data they are easy to understand and interpret which makes the
nonnegative factorization invaluable for scientific applications
\cite{Cichocki:Zdunek:Phan:Amari:2009}.
Classical tensor decompositions corresponding to nonnegative tensor
ranks are nonnegative Canonical Polyadic Decomposition
(nnCPD)~\cite{Cichocki:Zdunek:Phan:Amari:2009} and nonnegaive Tucker
Decomposition (nnTD)~\cite{Cichocki:Zdunek:Phan:Amari:2009}.
In Tucker, the minimum dimensions of the core tensor are often called
\emph{multirank} and the concept of \emph{nonnegative multirank} in
nnTD is introduced in Section~\ref{sec:Nonnegative_decompositions}.

In addition to nonnegativity, various other constraints on the
decomposition are often needed to take into account the available
domain knowledge and subject matter expertise and extract explainable
and meaningful latent features. The canonical decomposition with
linear constraints (CANDELINC) \cite{carroll1980candelinc} is one of
these decompositions. A preprocessing step in CANDELINC is Tucker
compression, which often leads to a larger residual error but also to
interpretable latent features \cite{bro1998improving}.

In this work, we derive formulation of nonnegative CANDELINC
(nnCANDELINC). This is accomplished via nnTDs we refer to as
\emph{minimal nonnegative Tucker Decompositions}: A minimal TD/nnTD is
where the Tucker core has the smallest shape possible (Definitions
\ref{def:minTD}, \ref{def:min_nTD}). In Section~\ref{sec:background},
we discuss the well-known fact that for real valued tensors minimal
Tuckers always exist \cite{Hackbusch:2012}, and preserve the rank of
the original tensor to the Tucker core (Theorem
\ref{thm:minTDExists}). We then relate the CPD to the minimal Tucker,
which leads to CANDELINC (Theorem \ref{def:cpd-td}). The previous
Theorem guarantees CANDELINC will successfully find a rank
factorization of the tensor. The nonnegative counterpart to CANDELINC
faces greater challenges however. In
Section~\ref{sec:Nonnegative_decompositions}, we discuss nnCPD with
rank deficiency (aka PARALIND for real valued tensors
\cite{Bro:Harshman:Sidiropoulos:Lundy:2009,wei2019compressed}), and
its relation to the minimal nnTD. We show that a minimal nnTD need not
exist (Example \ref{ex:cones_dne}) and even if it does exist, it need
not preserve the rank to the core (Example \ref{ex:equal}). However
under some mild conditions, some minimal nnTD will preserve the rank
(Theorem \ref{thm:minnTDExists}). Unfortunately, these conditions do
not guarantee that every minimal nnTD will preserve the rank (Example
\ref{ex:not_every_min_preserves_rank}). This naturally leads to the
discussion of when the nonnegative rank is preserved (Theorem
\ref{thm:nCPDmultirank}) and when we can overcome the challenges just
discussed. We therefore loosen our requirement on the shape of the
Tucker core. This leads us naturally to the definition of a canonical
nnTD (Definition \ref{def: canonical tucker}). We show that every
nonnegative tensor has a canonical nnTD which preserves the rank to
the core (Theorem \ref{thm: canonical preserve rank}). Finally, in
Section~\ref{sec:numerical}, we perform numerical experiments with
nnCANDELINC on synthetic and real-world datasets.  We consider two
different algorithms for nnCANDELINC: (i) Performing, first nnTD
compression, and then nnCPD on the core, and (ii) First nnCPD, and
then reconstruction of the linear dependence of the extracted factors
by Nonnegative Matrix Factorization (NMF). We also investigate the
effect of choosing the nonnegative canonical vs. nonnegative minimal
multirank.

\section{Decompositions of Real Valued Tensors}
\label{sec:background}
In this section, we review some of the basics of real-valued tensors
decompositions.
For notational simplicity, we consider only 3-way tensors, although
the analysis is valid for $d$-way tensors.
A detailed presentation of the basic results can be found
in~\cite{Cichocki:Zdunek:Phan:Amari:2009,Hackbusch:2012,Cohen:2016:PhD,Kolda-Bader:2009}. We
will begin with a few formal definitions on tensors and tensor
decompositions. This will include the concept of a minimal subspace,
which will motivate our definition of a minimal Tucker
decomposition. We will then go on to provide some results on rank
perseverance to the core of a minimal Tucker decomposition. Precise
notation for $n-$mode multiplication and unfolding used throughout the
text can be found in the Appendix.

\begin{definition}
  For vectors $\va^{(1)} \in \REAL^{N_1}$, $\va^{(2)} \in
  \REAL^{N_2}$, $\va^{(3)} \in \REAL^{N_3}$, the tensor product is the
  3-way tensor $\va^{(1)}\otimes\va^{(2)}\otimes\va^{(3)}$ given by
  \begin{align*}
    \left(\va^{(1)}\otimes\va^{(2)}\otimes\va^{(3)} \right)_{i,j,k} = \va^{(1)}_i \va^{(2)}_j \va^{(3)}_k.
  \end{align*}
  The tensor $\va^{(1)}\otimes\va^{(2)}\otimes\va^{(3)}$ is referred
  to as a \emph{rank-1}, elementary, or decomposable tensor. For $U_i$
  subspace of $\REAL^{N_i}$, the \emph{tensor product space} $U_1
  \otimes U_2 \otimes U_3$ consists of all linear combinations of
  elementary tensors where $\va^{(i)} \in U_i$.
\end{definition}
The tensor product space
$\REAL^{N_1}\otimes\REAL^{N_2}\otimes\REAL^{N_3}$ is isomorphic to the
linear space of 3-way arrays $\REAL^{N_1\times N_2\times N_3}$.
Thus for ease of notation, we will often write $\tX\in\REAL^{N_1\times
  N_2\times N_3}$ for a real 3-way tensor of dimension $N_1\times
N_2\times N_3$, with components $\tX=(\tX_{i,j,k})$, for $i,\,j,\,k$
ranging from $1$ to $N_1$, $N_2$, and $N_3$, respectively.
Every $\tX\in\REAL^{N_1}\otimes\REAL^{N_2}\otimes\REAL^{N_3}$ can be
written as
$\tX=\sum_{i,j,k}\tX_{i,j,k}\ve^{(1)}_{i}\otimes\ve^{(2)}_{j}\otimes^{(3)}_{k}$,
where, $\{\ve^{(1)}_i\}$, $\{\ve^{(2)}_j\}$, and $\{\ve^{(3)}_k\}$ are
the canonical basis vectors of $\REAL^{N_1}$, $\REAL^{N_2}$, and
$\REAL^{N_3}$, respectively.
However, every tensor $\tX \in \REAL^{N_1\times N_2\times N_3}$ can be
decomposed in many different ways.  And perhaps most significant is
the decomposition as a weighted sum of rank-1 tensors:

\begin{definition}
  For every tensor $\tX \in \REAL^{N_1\times N_2\times N_3}$, there
  exists a sufficiently large positive integer $r$ such that $\tX$ may
  be written as
  \begin{align}
    \tX=\sum_{n=1}^{r}\lambda_n\va^{(1)}_{n}\ccirc\va^{(2)}_{n}\ccirc\va^{(3)}_{n},
    \label{eq:pd}
  \end{align}
  where $\lambda_n\in\REAL$ and $\va^{(i)} \in \REAL^{N_i}$ are unit
  vectors.  Such a decomposition is a \emph{polyadic decomposition}.
  The \emph{rank} of a tensor is defined as the smallest integer
  number $r$ of rank-1 terms for which a polyadic decomposition
  exists, or
  \begin{align}
    \rank{(\tX)} = \min
    \bigg\{
    r~\big\vert~\tX = \sum_{n=1}^r\lambda_n\va^{(1)}_{n}\ccirc\va^{(2)}_{n}\ccirc\va^{(3)}_{n},
    \lambda_n\in\REAL,\,\va^{(i)}_n\in\REAL^{N_i},\,i=1,2,3
    \bigg\}.
  \end{align}
  A corresponding decomposition is called a \emph{Canonical Polyadic
  Decomposition (CPD) of $\tX$}.
\end{definition}

Collecting the vectors $\va^{(i)}_n$ into factor matrices
$\mA^{(i)}=\big[\;\va^{(i)}_1\;\big\vert\;\hdots\;\big\vert\;\va^{(i)}_r\;\big]$
and the coefficients $\lambda_n$ into a superdiagonal tensor $\tD$
allows us to represent CPD as the product of a superdiagonal tensor
$\tD$ and factor matrices, or
\begin{align}
  \tX=\tD\times_1\mA^{(1)}\times_2\mA^{(2)}\times_3\mA^{(3)},
  \label{eq:CPD:matrix-format}
\end{align}
as seen in Figure \ref{fig:nCPD}, panel A.  Here the $n-$mode
multiplication, $\times_n$ is defined in the Appendix ( Definition
\ref{def:n-mode}).  In general, $\tX$ does not require the full
ambient space $\REAL^{N_1}\otimes\REAL^{N_2}\otimes\REAL^{N_3}$ to
represent it.
Indeed, $\tX$ can be contained in the tensor product of subspaces
$\sU_1\otimes\sU_2\otimes\sU_3$ where $\sU_i$ is a subspace of
$\REAL^{N_i}$.  This is the concept behind a Tucker Decomposition:

\begin{figure}[!t]
  \centering
  \includegraphics[width=0.8\textwidth]{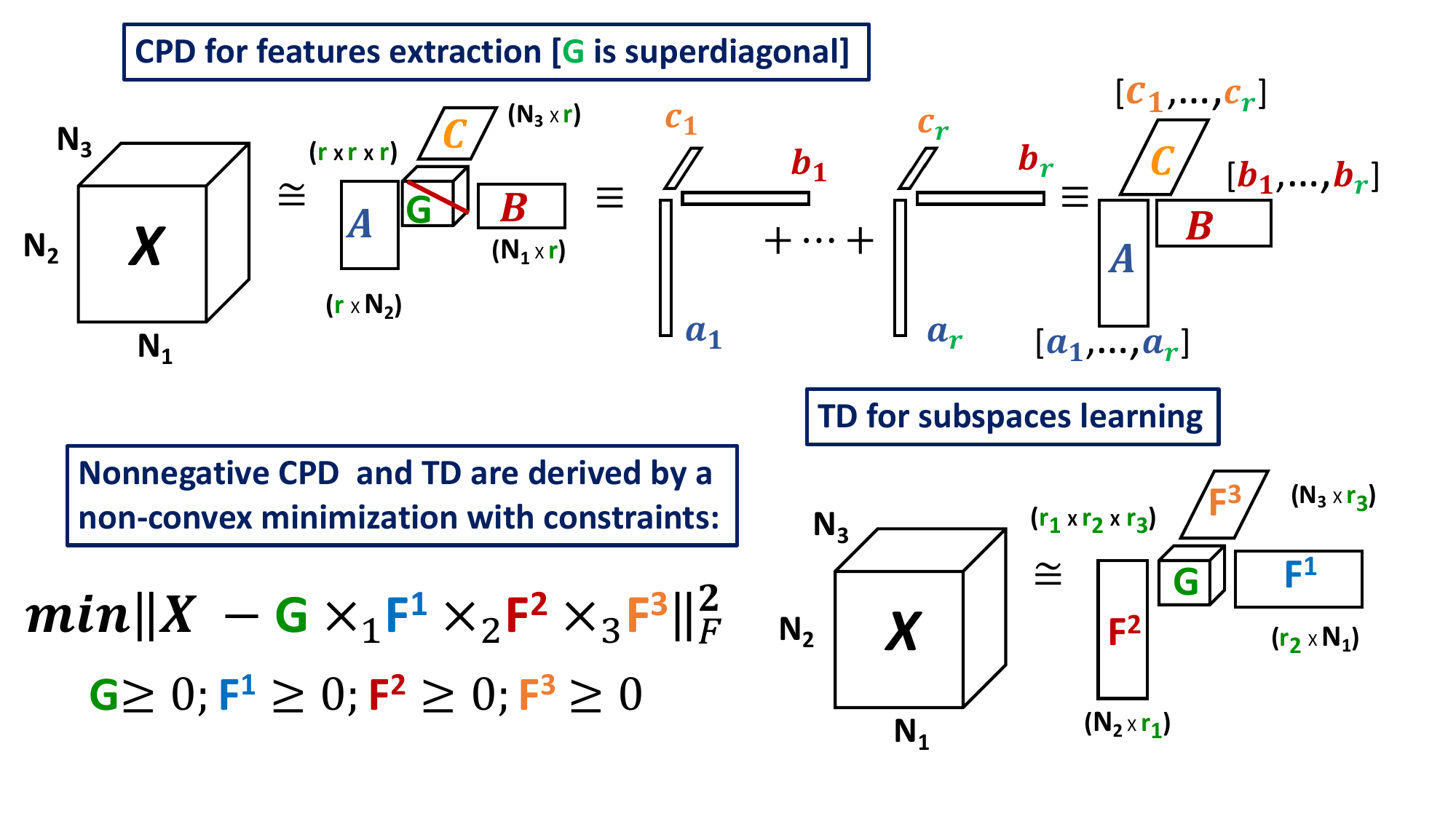}
  \caption{Two classical tensor decompositions: A) Canonical Polyadic
    Decomposition (CPD) of a 3-dimensional tensor $\tX$ of size
    $N_1\times N_2\times N_3$ into a superdiagonal core tensor $\tG
    \equiv \tD$ of size $r \times r \times r$ and three matrix
    factors, $\mA$, $\mB$, and $\mC$. B) Tucker Decomposition (TD) of
    a 3-dimensional tensor $\tX$ into a dense core tensor $\tG$ of
    size $r_1 \times r_2 \times r_3$ and three matrix factors,
    $\mF^{(1)}$, $\mF^{(2)}$, and $\mF^{(3)}$.}
  \label{fig:nCPD}
\end{figure}

\begin{definition}
  The \emph{Tucker Decomposition} (TD) is a weighted tensor product
  decomposition of the form
  \begin{align}
    \tX=\sum^{r_1,r_2,r_3}_{n_1,n_2,n_3=1}\tG_{n_1,n_2,n_3}\vf^{(1)}_{n_1}\ccirc\vf^{(2)}_{n_2}\ccirc\vf^{(3)}_{n_3},
    \label{eq:TD:componentwise}
  \end{align}
  where the vectors $\vf^{(i)}\in\REAL^{N_i}$, for $i=1,2,3$, and the
  core tensor $\tG\in\REAL^{r_1}\otimes\REAL^{r_2}\otimes\REAL^{r_3}$.
\end{definition}
Tucker decomposition factorizes tensor $\tX$ into the product of a
tensor core $\tG$ and three factor matrices
$\mF^{(i)}=\begin{bmatrix}\vf^{(i)}_1\vert\hdots\vert\vf^{(i)}_{r_i}\end{bmatrix}\in\REAL^{N_i\times
  r_i}$, for $i=1,2,3$, as seen in Figure \ref{fig:nCPD} panel B.
Similarly to~\eqref{eq:CPD:matrix-format}, we can
reformulate~\eqref{eq:TD:componentwise} as
\begin{align}
  \tX=\tG\times_1\mF^{(1)}\times_2\mF^{(2)}\times_3\mF^{(3)}\;.
  \label{eq:TD:matrix-format}
\end{align}
For a tensor $\tX\in\REAL^{N_1}\times \REAL^{N_2}\times\REAL^{N_3}$,
the matrix factors $\mF^{(i)}$ in a Tucker decomposition are
associated with such subspaces of $\REAL^{N_i}$.  Given a matrix
$\mF\in\REAL^{N\times r}$ we let $\col(F)$ denote the column space of
$F$.  Then the following is a direct consequence of Equation
\ref{eq:TD:componentwise}:

\begin{proposition}\label{prop:tdspaces}
  Given three matrices $\mF^{(1)}$, $\mF^{(2)}$, and $\mF^{(3)}$, a
  tensor $\tX$ admits the Tucker decomposition
  $\tX=\tG\times_1\mF^{(1)}\times_2\mF^{(2)}\times_3\mF^{(3)}$ if and
  only if
  $\tX\in\col(\mF^{(1)})\otimes\col(\mF^{(2)})\otimes\col(\mF^{(3)})$.
\end{proposition}

Every Tucker decomposition of a three-way tensor is linked to three
integer numbers, namely, $r_1$, $r_2$, and $r_3$, from
$\tG\in\REAL^{r_1\times r_2\times r_3}$.
The immediate question is: what are the permissible or minimal values
of $r_1$, $r_2$, and $r_3$ such that there exists a Tucker
decomposition with a core tensor of these dimensions?
The smallest such values would describe the maximal permissible
lossless compression within the shape of the tensor.  This information
is encoded in the concept of the \emph{minimal subspaces} and
\emph{minimal tensor multirank}\cite{Hackbusch:2012}.

\begin{definition}
  \label{def:minimal:subspaces}
  Given a tensor
  $\tX\in\REAL^{N_1}\otimes\REAL^{N_2}\otimes\REAL^{N_3}$, the
  \emph{minimal subspaces associated with $\tX$} are subspaces
  $\sU_i^{\min}\subset\REAL^{N_i}$ such that
  $\tX\in\sU_1^{\min}\otimes\sU_2^{\min}\otimes\sU_3^{\min}$ and if
  $\tX\in\sU_1\otimes\sU_2\otimes\sU_3$ then
  $\sU_i^{\min}\subset\sU_i$.
\end{definition}

We remark that minimal subspaces always exist and are unique.  
Indeed, one can show that
\begin{equation}
  \label{equation:intersection}
  \left(\sU_1\otimes\sU_2\otimes\sU_3\right)
  \bigcap
  \left(\sU_1'\otimes\sU_2'\otimes\sU_3'\right) 
  = \bigotimes_{i=1}^3\,\sU_i\cap\sU_i'
\end{equation}
for any collection of subspaces $\sU_i$, $\sU_i'\subset\REAL^{N_i}$
\cite{Hackbusch:2012}.
It then follows that
\begin{align*}
  \sU_1^{\min}\otimes\sU_2^{\min}\otimes\sU_3^{\min} = \bigcap \left\{ \sU_1\otimes\sU_2\otimes\sU_3 : \tX \in \sU_1\otimes\sU_2\otimes\sU_3  \right\}.
\end{align*}
Eq.~\eqref{equation:intersection} also shows that the minimal
subspaces can be found coordinatewise rather than simultaneously.
Hence,  if $\sU_i^{\min}$ are the minimal subspaces found such that
\begin{align*}
  \tX \in \sU_1^{\min}\otimes\REAL^{N_2}\otimes\REAL^{N_3}; \quad
  \tX \in \REAL^{N_1}\otimes\sU_2^{\min}\otimes\REAL^{N_3}; \quad
  \tX \in \REAL^{N_1}\otimes\REAL^{N_2}\otimes\sU_3^{\min},
\end{align*}
then by Eq.~\eqref{equation:intersection}, we have
$\tX\in\sU_1^{\min}\otimes\sU_2^{\min}\otimes\sU_3^{\min}$. Associated
with the minimal subspaces of $\tX$ is the concept of the \emph{$i$-th
minimal multirank} of $\tX$.
\begin{definition}
  The \emph{$i$-th minimal multirank} of a tensor $\tX$, denoted by
  $\multirank_i(\tX)$, is the dimension of the $i$-th minimal subspace
  $\sU_i^{\min}$.
  The \emph{minimal multilinear rank} of $\tX$ is the triple of dimensions
  \begin{align*}
    \multirank(\tX) = \big(\multirank_1(\tX),\,\multirank_2(\tX),\,\multirank_3(\tX)\big).
  \end{align*}
\end{definition}

We note that the $i$-th minimal multirank of $\tX$ does not depend on
the $j$-th tensor coordinate for $j\neq i$.
Formally, the first minimal multirank of $\tX$ is given by
\begin{align*}
  \multirank_1(\tX) = 
  \min\big\{
  \operatorname{dim}(\sU_1)~\vert~\tX\in\sU_1\otimes\REAL^{N_2}\otimes\REAL^{N_3},\sU_1\subset\REAL^{N_1}
  \big\}
\end{align*}
with analogous definitions for the second and the third minimal
multiranks.
For any Tucker decomposition
$\tX=\tG\times_1\mF^{(1)}\times_2\mF^{(2)}\times_3\mF^{(3)}$, it holds
that $\rank(\mF^{(i)})\geq\multirank_{i}(\tX)$, $i=1,2,3$, since by
Proposition \ref{prop:tdspaces} the span of the columns of matrix
factor $\mF^{(i)}$ must contain the corresponding minimal subspace
$\sU_i^{\min}$, i.e.,
$\sU_i^{\min}\subset\operatorname{span}(\mF^{(i)})$.
We record the following well known connection between the minimal
multirank of $\tX$ and its unfoldings (for definition of unfolding see
Appendix, Definition \ref{def:unfolding}).

\begin{proposition}
\label{prop:unfold_multirank}
  Given a tensor $\tX$, $\sU_i^{\min} = \col(\matricization_i(\tX))$
  and $\multirank_{i}(\tX) = \rank(\matricization_i(\tX))$
  \cite{Hackbusch:2012}.
\end{proposition}

In general, a Tucker decomposition does not satisfy the identity
$\rank(\mF^{(i)})=\multirank_{i}(\tX)$.
This fact motivates us to introduce the notion of \emph{ minimal}
  Tucker Decomposition (\emph{minimal} TDs) in the next definition, which is
a Tucker decomposition with core dimensions corresponding to the
\emph{minimal multirank}.
\begin{definition}\label{def:minTD}
  Consider tensor $\tX\in\REAL^{N_1\times N_2\times N_3}$.
  We say that the Tucker decomposition
  $\tX=\tG\times_1\mF^{(1)}\times_2\mF^{(2)}\times_3\mF^{(3)}$ is
  \emph{minimal} if the dimensions of the core tensor $\tG$ are equal
  to the minimal multiranks, i.e.,
  $\tG\in\REAL^{\multirank_1(\tX)\times\multirank_2(\tX)\times\multirank_3(\tX)}$
  and $\mF^{(i)}\in\REAL^{N_i\times\multirank_i(\tX)}$, $i=1,2,3$.
\end{definition}
In a minimal TD, $\mF^{(i)}\in\REAL^{N_i\times\multirank_i(\tX)}$
implies that $\rank(\mF^{(i)}) \leq \multirank_i(\tX)$.  However by
the discussion above, $\rank(\mF^{(i)}) \geq \multirank_i(\tX)$, so
that $\rank(\mF^{(i)}) = \multirank_i(\tX)$. This is rather different
than the case of the loading matrices in a CPD, where rank deficiency
can occur. The following simple example demonstrates that a CPD need
not be minimal TD.

\begin{example}
  \label{ex:real_minTD_not_cpd}
  Let $\tX \in \REAL^{2,2,2}$ be the rank 2 tensor
  \[
  \tX = \left[
    \begin{array}{cc|cc}
      0 &\   1 &\  1 &\  0\\
      0 &\   1 &\  1 &\ 0
    \end{array}
    \right].
  \]
  It is not hard to see that $\multirank_1(\tX) = 1$, while
  $\multirank_2(\tX) = \multirank_3(\tX) =2$.  Hence any minimal TD will
  satisfy $\tG \in \REAL^{1,2,2}$.  However since $\tX$ is rank 2, the
  rank decomposition will have the shape $\tD \in \REAL^{2,2,2}$. Thus,
  a CPD of $\tX$ does not need to be a minimal Tucker decomposition.
\end{example}

\subsection{Real Rank Preservation to Minimal Tucker Core}
Given a Tucker decomposition $\tX
=\tG\times_1\mF^{(1)}\times_2\mF^{(2)}\times_3\mF^{(3)}$, we say that
the rank is preserved to the core if $\rank(\tX) = \rank(\tG)$.  Not
every TD needs to preserve the rank to the Tucker core, as the
following simple example illustrates:

\begin{example}
  Let $\tX \in \REAL^{2,2,2}$ be the rank 1 tensor
  \[
  \tX = \left[
    \begin{array}{cc|cc}
      1 &\   1 &\  1 &\  1\\
      1 &\   1 &\  1 &\ 1
    \end{array}
    \right].
  \]
  Then $\tX$ can be decomposed using a rank $2$ core as
  \[
  \tX  = \left[
    \begin{array}{cc|cc}
      1 &\   0 &\  1 &\  0\\
      0 &\   1 &\  0 &\ 1
    \end{array}
    \right]
  \times_1
  \begin{bmatrix}
    1 & 1 \\
    1 & 1
  \end{bmatrix}.
  \]
\end{example}

One always has the following simple rank and minimal Tucker
relationships which we record in the following lemma:

\begin{lemma}
  \label{lemma:rank and min exists}
  For any Tucker decomposition $\tX
  =\tG\times_1\mF^{(1)}\times_2\mF^{(2)}\times_3\mF^{(3)}$, $\rank(\tX)
  \leq \rank(\tG)$.  Moreover, a minimal Tucker decomposition always
  exists for real factorizations.
\end{lemma}
\begin{proof}
  Suppose $\tX =\tG\times_1\mF^{(1)}\times_2\mF^{(2)}\times_3\mF^{(3)}$
  is a Tucker decomposition of $\tX$.  Then consider a CPD of $\tG$
  given by
  \begin{align}
    \tG=\tD_{\tG}\times_1\mA_{\tG}^{(1)}\times_2\mA_{\tG}^{(2)}\times_3\mA_{\tG}^{(3)}.
    \label{eq:minTDExists:G-CPD}
  \end{align}
  Then, substituting the CPD of $\tG$ into the Tucker decomposition of $\tX$ yields
  \begin{align*}
    \begin{array}{rll}
      \tX 
      &= \tG\times_1\mF^{(1)}\times_2\mF^{(2)}\times_3\mF^{(3)}
      \nonumber\\[0.5em]
      &= 
      \tD_{\tG}
      \times_1\mA_{\tG}^{(1)}
      \times_2\mA_{\tG}^{(2)}
      \times_3\mA_{\tG}^{(3)}
      \times_1\mF^{(1)}
      \times_2\mF^{(2)}
      \times_3\mF^{(3)}
      \nonumber\\[0.5em]
      &= 
      \tD_{\tG}
      \times_1\big(\mF^{(1)}\mA_{\tG}^{(1)}\big)
      \times_2\big(\mF^{(2)}\mA_{\tG}^{(2)}\big)
      \times_3\big(\mF^{(3)}\mA_{\tG}^{(3)}\big).
    \end{array}
  \end{align*}
  The last right-hand side is a polyadic decomposition of $\tX$ with
  $\rank(\tG)$ summands proving that $\rank(\tX)\leq\rank(\tG)$.

  It is also not hard to see that a minimal Tucker always exists for
  real factorizations.  Let $\mF^{(i)}$ be basis matrices for
  $\sU_i^{\min}$.
  Then by definition,
  \[
  \tX \in \sU_1^{\min}\otimes\sU_2^{\min}\otimes\sU_3^{\min} = \col(\mF^{(1)})\otimes\col(\mF^{(2)})\otimes\col(\mF^{(3)}).
  \]
  By Proposition \ref{prop:tdspaces}, there exists a $\tG$ such that
  $\tX =\tG\times_1\mF^{(1)}\times_2\mF^{(2)}\times_3\mF^{(3)}$.  As
  $\mF^{(i)}\in\REAL^{\multirank_i(\tX) \times \sN_i}$, this is a
  minimal TD.
\end{proof}

Lemma \ref{lemma:rank and min exists} demonstrates that a minimal
Tucker decomposition can be constructed by choosing basis matricies
for $\sU_i^{\min}$. However, this is how all real minimal Tucker
decompositions are formed.  Indeed if
$\tX=\tG\times_1\mF^{(1)}\times_2\mF^{(2)}\times_3\mF^{(3)}$ is a
minimal Tucker, then by Proposition \ref{prop:tdspaces},
$\tX\in\col(\mF^{(1)})\otimes\col(\mF^{(2)})\otimes\col(\mF^{(3)})$.
Since $\rank(\mF^{(i)}) = \multirank_i(\tX)$ and
$\mF^{(i)}\in\REAL^{\multirank_i(\tX) \times \sN_i}$, this implies
$\mF^{(i)}$ are basis matrices for $\sU_i^{\min}$.

While a minimal Tucker decomposition always exists, what we are
interested in is the preservation of the rank of $\tX$ to the minimal
core $\tG$.
The next theorem establishes that minimal Tucker decompositions do
always preserve the rank to the core. While we believe this result is
known, it does not appear to be written explicitly down in the
literature. Hence, we record it alongside its proof:
\begin{theorem}\label{thm:minTDExists}
  Given a real tensor $\tX$, any minimal Tucker decomposition
  $\tX=\tG\times_1\mF^{(1)}\times_2\mF^{(2)}\times_3\mF^{(3)}$
  satisfies $\rank(\tX)=\rank(\tG)$.
\end{theorem}
\begin{proof}
  The proof is constructive.
  Let $\sU_i^{\min}$ be the $i$-th minimal subspace of $\tX$ according
  to the definition of the minimal multirank.
  For every $i=1,2,3$, we choose the basis matrix
  $\mF^{(i)}\in\REAL^{\multirank_i(\tX) \times \sN_i}$ such that
  $\operatorname{span}(\mF^{(i)})=\sU_i^{\min}$.
  By construction, $\mF^{(i)}$ is a full column rank matrix, its rank
  being equal to $\multirank_i(\tX)$. 
  Hence, $\mF^{(i)}$ admits the (left) Moore-Penrose pseudo-inverse
  ${\mF^{(i)}}^{\dagger}\in\REAL^{\sN_i \times \multirank_i(\tX)}$, so
  that ${\mF^{(i)}}^{\dagger}{\mF^{(i)}}=\mI$.
  Moreover, the matrix $\mP^{(i)}:= \mF^{(i)} {\mF^{(i)}}^{\dagger}$
  is the projection onto the column space of $\mF^{(i)}$ (which is
  $\sU_i^{\min}$).
  Let 
  \begin{align}
    \tG
    = \tX\times_1{\mF^{(1)}}^{\dagger}\times_2{\mF^{(2)}}^{\dagger}\times_3{\mF^{(3)}}^{\dagger}.
    \label{eq:minTDExists:G-def}
  \end{align}
  We show that
  $\tX=\tG\times_1\mF^{(1)}\times_2\mF^{(2)}\times_3\mF^{(3)}$.
  Indeed, by redistributing we verify that
  \begin{align*}
    \tG\times_1\mF^{(1)}\times_2\mF^{(2)}\times_3\mF^{(3)} 
    = 
    \tX
    \times_1\mP^{(1)} 
    \times_2 \mP^{(2)}
    \times_3\mP^{(3)}.
  \end{align*}
  Hence, it suffices to show that $\tX \times_i \mP^{(i)} = \tX$ for
  each $i$. 
  This happens if and only if
  $\col(\matricization_i(\tX))\subset\col(\mP^{(i)})=\col(\mF^{(i)})=\sU_i^{\min}$.
  By Proposition \ref{prop:unfold_multirank}, these are equal.
  Therefore,
  $\tG\times_1\mF^{(1)}\times_2\mF^{(2)}\times_3\mF^{(3)}=\tX$.
  
  To prove the identity $\rank(\tX)=\rank(\tG)$, it suffices to show
  that $\rank(\tX) \geq \rank(\tG)$.
  Let
  \begin{align}
    \tX = 
    \tD_{\tX}
    \times_1\mA_{\tX}^{(1)}
    \times_2\mA_{\tX}^{(2)}
    \times_3\mA_{\tX}^{(3)}
    \label{eq:minTDExists:X-CPD}
  \end{align}
  be a CPD of $\tX$ with matrix factors
  $\mA_{\tX}^{(i)}\in\REAL^{N_i\times r}$ and superdiagonal core
  tensor $\tD_{\tX}\in\REAL^{r\times r\times r}$, where
  $r=\rank(\tX)$.
  Then, starting from~\eqref{eq:minTDExists:G-def} and
  using~\eqref{eq:minTDExists:X-CPD}, a straightforward calculation
  yields:
  \begin{align*}
    \begin{array}{rll}
      \tG 
      &= \tX
      \times_1{\mF^{(1)}}^{\dagger}
      \times_2{\mF^{(2)}}^{\dagger}
      \times_3{\mF^{(3)}}^{\dagger}
      \nonumber\\[0.5em]
      &= \tD_{\tX}
      \times_1\mA_{\tX}^{(1)}
      \times_2\mA_{\tX}^{(2)}
      \times_3\mA_{\tX}^{(3)}
      \times_1\big({\mF^{(1)}}^{\dagger}\big)
      \times_2\big({\mF^{(2)}}^{\dagger}\big)
      \times_3\big({\mF^{(3)}}^{\dagger}\big)
      \nonumber\\[0.5em]
      &= 
      \tD_{\tX}
      \times_1\big({\mF^{(1)}}^{\dagger}\mA_{\tX}^{(1)}\big)
      \times_2\big({\mF^{(2)}}^{\dagger}\mA_{\tX}^{(2)}\big)
      \times_3\big({\mF^{(3)}}^{\dagger}\mA_{\tX}^{(3)}\big).
    \end{array}
  \end{align*}
  The last right-hand side is a polyadic decomposition of $\tG$ with
  $\rank(\tX)$ summands, proving that $\rank(\tG)\leq\rank(\tX)$.
\end{proof}

\subsection{Rank Deficiency in CPD Factors and the CANDELINC Solution}
In this subsection, we discuss the challenges that arise from rank
deficient factors in a CPD and how the CANDELINC method
~\cite{Bro:Harshman:Sidiropoulos:Lundy:2009, Kolda-Bader:2009} can
provide a suitable decomposition. Theorem \ref{thm:minTDExists}
establishes that one can always construct minimal Tucker
decompositions that will preserve the rank.  The following theorem
relates the uniqueness of a CPD to the minimal TDs:

\begin{theorem}[Ranks of CPD factors related to minimal TD]
  \label{def:cpd-td}
  Let $\tX=\tD\times_1\mA^{(1)}\times_2\mA^{(2)}\times_3\mA^{(3)}$ be
  a CPD of $\tX$, then $\multirank_{i}(\tX) \leq
    \rank(\mA^{(i)})$.
  Furthermore, if the CPD is unique, then
  $\multirank_{i}(\tX)=\rank(\mA^{(i)})$.  In this case, $\col(\mA^{(i)}) = \sU_i^{\min}$.
\end{theorem}
\begin{proof}
  Let $\tX=\tD\times_1\mA^{(1)}\times_2\mA^{(2)}\times_3\mA^{(3)}$
  be a CPD of $\tX$.  Since
  $\tD\times_1\mA^{(1)}\times_2\mA^{(2)}\times_3\mA^{(3)}$ is a
  Tucker decomposition, by definition of minimal multirank we have
  $\multirank_{i}(\tX) \leq \rank(\mA^{(i)})$. Now suppose that
  $\tX$ has a unique CPD, and consider a minimal TD of the form
  $\tX=\tG\times_1\mF^{(1)}\times_2\mF^{(2)}\times_3\mF^{(3)}$.
  We recall that $\rank(\mF^{(i)})=\multirank_{i}(\tX)$ by definition.
  On its turn, the Tucker core admits the CPD
  $\tG=\tD_{\tG}\times_1\mB^{(1)}\times_2\mB^{(2)}\times_3\mB^{(3)}$.
  Substituting the CPD of $\tG$ in the TD of $\tX$, we obtain the
  alternative CPD
  $\tX=\tD_{\tG}\times_1\mF^{(1)}\mB^{(1)}\times_2\mF^{(2)}\mB^{(2)}\times_3\mF^{(3)}\mB^{(3)}$.
  Since we assume that the CPD is unique, with appropriate scalings
  and permutations, which are rank-preserving operations, we obtain
  that $\mA^{(i)}=\mF^{(i)}\mB^{(i)}$.
  Therefore,
  $\rank(\mA^{(i)})\leq\rank(\mF^{(i)})=\multirank_{i}(\tX)$. It
  follows from Proposition \ref{prop:unfold_multirank} and rank
  arguments that $\sU_i^{\min} = \col(\matricization_i(\tX)) =
  \col(\mA^{(i)})$.
\end{proof}

Theorem~\ref{def:cpd-td} suggests why a direct CPD computation can be
algorithmically problematic.
Let $\tX=\tD\times_1\mA^{(1)}\times_2\mA^{(2)}\times_3\mA^{(3)}$ be
the unique CPD of a rank $r$ tensor with $\multirank_i(\tX) = r_i$.
If $r_i<r$, as is the case with probability 1 for many shaped tensors
\cite{qi2018very,Qi-Comon-Lim:2016}, then $\mA^{(i)}$ is a rank
deficient matrix by Theorem~\ref{def:cpd-td}.
Indeed, $\mA^{(i)}$ is an ($N_i\times r$)-sized matrix with only $r_i$
linearly independent columns.
Algorithmically, finding rank deficient matrices without an explicit
rank constraint for tensors of the form $\tX=\widetilde{\tX}+\tE$ is
challenging, as the rank deficient subspaces of the factors of
$\widetilde{\tX}$ can always be expanded to accommodate some of the
noise, $\tE$.

The proof of Theorem 2 suggests a more suitable method for computing
the CPD of $\tX$. First compute a minimal TD of $\tX$ (which will
preserve the rank); then compute a CPD of the TD core (which will lack
rank deficiency); and finally substitute the CPD of the TD core into
the TD and obtain a CPD of the original tensor.
Bro \textit{et al.} followed this strategy in their construction of
the PARALIND models, cf.~\cite{Bro:Harshman:Sidiropoulos:Lundy:2009},
and Carroll \textit{et al.} followed this strategy in their
construction of the CANDELINC models, cf.~\cite{carroll1980candelinc}.
Formally, if
$\tX=\tG\times_1\mF^{(1)}\times_2\mF^{(2)}\times_3\mF^{(3)}$ is a
minimal TD, and
$\tG=\tD_{\tG}\times_1\mA_{\tG}^{(1)}\times_2\mA_{\tG}^{(2)}\times_3\mA_{\tG}^{(3)}$
is the CPD of the Tucker core,
then each factor $\mA_{\tG}^{(i)}$ is a full column rank matrix, avoiding
the algorithmic problems previously discussed.
A simple substitution yields a CPD of $\tX$ where each loading matrix
is rank factored, i.e., $\tX = \tD_{\tG}
\times_1\big(\mF^{(1)}\mA_{\tG}^{(1)}\big)
\times_2\big(\mF^{(2)}\mA_{\tG}^{(2)}\big)
\times_3\big(\mF^{(3)}\mA_{\tG}^{(3)}\big)$, and we have explicitly
the linear constraints of the CPD factors.

\section{Nonnegative Decompositions of Nonnegative Tensors}
\label{sec:Nonnegative_decompositions}
Following~ \cite{landsberg2012tensors,comon2014tensors},
we now present the nonnegative counterparts to the discussion for real
tensors above.  This theory necessarily depends on some basic
knowledge of nonnegative matrix factorizations. For the unfamiliar
reader, we have provided some background information in the Appendix.
Throughout, we let $\REAL_+$ denote the nonnegative real numbers.  All
of the basic definitions from real tensors will carry over to
nonnegative with some appropriate adaptations.  While real rank
factorizations fundamentally rely on subspaces, nonnegative
factorizations are concerned with the nonnegative analog of subspaces
- polyhedral cones.

\begin{definition}
  A \textit{convex cone} is a subset $C \subset \REAL_+^N$ that is
  closed under addition of vectors and $\REAL_+$ scalar multiplication.
  Given $W \subset \REAL_+^N$, the non-negative span of $W$ defines a
  cone.  A subset of the cone $W \subset C \subset \REAL_+^N$ is a
  \textit{generating set} if its span is equal to $C$. The
  \textit{order} of the cone $C \subset \REAL_+^N$, denoted
  $\mathcal{O}(C)$, is the size of a minimal generating set. A cone is
  \textit{polyhedral} if $\mathcal{O}(C) < \infty$. Given a nonnegative
  matrix $W \in \REAL_+^{N,R}$, we define the \textit{cone of the matrix
    $W$} to be
  \[
  \mbox{cone}(W) = \{W h: h \in \REAL_+^R\} \subset \REAL_+^N.
  \]
\end{definition}

Every polyhedral cone $C \subset \REAL_+^N$ is $\cone(W)$ for some
nonnegative matrix $W \in \REAL_+^{N,R}$. Furthermore, every
polyhedral cone can be equivalently described as the intersection of
half spaces \cite{Bertsimas-Tsitsiklis:1997}.  With the precise
definition of cone, we can now define the analogous tensor product
space of cones, and the associated nonnegative tensor decompositions.

\begin{definition}
  For vectors $\va^{(1)} \in \REAL_+^{N_1}$, $\va^{(2)} \in
  \REAL_+^{N_2}$, $\va^{(3)} \in \REAL_+^{N_3}$, the tensor product is
  the 3-way tensor $\va^{(1)}\otimes\va^{(2)}\otimes\va^{(3)}$ given
  by
  \begin{align*}
    \left(\va^{(1)}\otimes\va^{(2)}\otimes\va^{(3)} \right)_{i,j,k} =
    \va^{(1)}_i \va^{(2)}_j \va^{(3)}_k.
  \end{align*}
 The tensor $\va^{(1)}\otimes\va^{(2)}\otimes\va^{(3)}$ is referred to
 as a nonnegative \emph{rank-1}, elementary, or decomposable
 tensor. For $\cC_i$ a polyhedral cone of $\REAL_+^{N_i}$, the
 \emph{tensor product space} $\cC_1 \otimes \cC_2 \otimes \cC_3$
 consists of all nonnegative linear combinations of elementary tensors
 where $\va^{(i)} \in \cC_i$.
\end{definition}

Analogous to the real case, every tensor $\tX \in \REAL^{N_1 \times
  N_2 \times N_3}$ can be decomposed in different ways.  The
definitions of polyadic and Tucker decompositions for tensors will
translate with the appropriate nonnegative adjustments.

\begin{definition}
  For every tensor $\tX \in \REAL_+^{N_1 \times N_2 \times N_3}$,
  there exists a sufficiently large positive integer $r$ such that
  $\tX$ may be written as
  \begin{align}
    \tX=\sum_{n=1}^{r}\lambda_n\va^{(1)}_{n}\ccirc\va^{(2)}_{n}\ccirc\va^{(3)}_{n},
    \label{eq:npd}
  \end{align}
  where $\lambda_n\in\REAL_+$ and $\va^{(i)} \in \REAL_+^{N_i}$ are
  unit vectors.  Such a decomposition is a \emph{nonnegative polyadic
  decomposition}.  The \emph{nonnegative rank} of a tensor is defined
  as the smallest integer number $r$ of rank-1 terms for which a
  polyadic decomposition exists, or
  \begin{align}
    \rank{(\tX)} = \min
    \bigg\{
    r~\big\vert~\tX = \sum_{n=1}^r\lambda_n\va^{(1)}_{n}\ccirc\va^{(2)}_{n}\ccirc\va^{(3)}_{n},
    \lambda_n\in\REAL_+,\,\va^{(i)}_n\in\REAL_+^{N_i},\,i=1,2,3
    \bigg\}.
  \end{align}
  A corresponding decomposition is called a \emph{nonnegative
  Canonical Polyadic Decomposition (nnCPD) of $\tX$}. For brevity, if
  the nonnegative qualifier is clear from context we may omit it when
  discussing various nonnegative ranks.
\end{definition}

It is immediately clear that for tensors $\rank_+(\tX)\geq\rank(\tX)$,
as the nnCPD is also a polyadic decomposition.
Analogous to the real case, $\tX$ does not require the full ambient
space of $\REAL_+^{N_1} \otimes \REAL_+^{N_2} \otimes \REAL_+^{N_3}$
to represent it.  It is possible that $\tX$ can be contained in the
tensor product of cones $\cC_1 \otimes \cC_2 \otimes \cC_3$ where
$\cC_i$ is a polyhedral cone of $\REAL_+^{N_i}$.  This once again
motivates the concept of a nonnegative Tucker decomposition:

\begin{definition}
  A \emph{nonnegative Tucker decomposition} (nnTD) of a nonnegative
  tensor is a nonnegative weighted tensor product decomposition of the
  form,
  \begin{align}
    \tX = \sum^{r_1,r_2,r_3}_{n_1,n_2,n_3=1} \tG_{n_1,n_2,n_3}
    \vf^{(1)}_{n_1}\ccirc\vf^{(2)}_{n_2}\ccirc\vf^{(3)}_{n_3}\;,
  \end{align}
  where the vectors $\vf^{(i)}_{n_i} \in \REAL^{N_i}_+$, for
  $i=1,2,3$, and the core tensor $\tG_{n_1,n_2,n_3} \in \REAL^{r_1}_+
  \otimes \REAL^{r_2}_+ \otimes \REAL^{r_3}_+$.
\end{definition}
The factors of an nnTD are associated with nonnegative cones, and are
inherently tied to the tensor belonging to the tensor product space of
these cones:

\begin{proposition}\label{prop:ntdCones}
  Given nonnegative matrices $\mF^{(1)}$, $\mF^{(2)}$, and
  $\mF^{(3)}$, a tensor admits an nnTD: $\tX = \tG \times_1 \mF^{(1)}
  \times_2 \mF^{(2)} \times_3 \mF^{(3)} \textrm{~if~and~only~if~} \tX
  \in \cone(\mF^{(1)}) \otimes \cone(\mF^{(2)}) \otimes
  \cone(\mF^{(3)})$.
\end{proposition}

One important subtle difference between polyhedral cones and subspaces
is that cone intersection does not commute with the tensor product.
That is, if $\cC_i, \cC_i' \subset \REAL_+^{N_i}$ are cones for
$i=1,2,3$, then
\begin{align*}
  \left( \cC_1\otimes \cC_2\otimes \cC_3 \right) 
  \bigcap \left( \cC_1'\otimes \cC_2'\otimes \cC_3' \right) 
  \neq \bigotimes_{i=1}^3 \cC_i\cap \cC_i'.
\end{align*}
Example \ref{ex:not_in_cones} shows that we cannot simply take the
``smallest'' cones via intersection as we could with subspaces.

\begin{example}\label{ex:not_in_cones}
  Consider a $3\times 3\times 2$ nonnegative tensor with the
  unfoldings,
  \begin{align*}
    \matricization_1(\tX) = \left[
      \begin{array}{ccc|ccc}
        1 &\  1 &\  1 &\  0 &\  0 &\  1\\
        1 &\  1 &\  2 &\  1 &\  1 &\  1\\
        1 &\  1 &\  2 &\  1 &\  1 &\  1
      \end{array}
    \right]\;,\\
    \matricization_2(\tX) = \left[
      \begin{array}{ccc|ccc}
        1 &\  1 &\  1 &\  0 &\  1 &\  1\\
        1 &\  1 &\  1 &\  0 &\  1 &\  1\\
        1 &\  2 &\  2 &\  1 &\  1 &\  1
      \end{array}
    \right]\;,\\
    \matricization_3(\tX) = \left[
      \begin{array}{ccc|ccc|ccc}
        1 &\  1 &\  1 &\  1 &\  1 &\  1 &\  1 &\  2 &\  2\\
        0 &\  1 &\  1 &\  0 &\  1 &\  1 &\  1 &\  1 &\  1
      \end{array}
    \right]\;.
  \end{align*}
  One can easily verify that:
  \begin{align*}
    \tX \in \cC_1 \otimes \REAL^{3}_+ \otimes \REAL^{2}_+ 
    \quad\text{where}\quad
    & \cC_1 = \operatorname{cone}( W^{(1)} )
    \quad\text{and}\quad
    W^{(1)} = \left[ \begin{array}{cc}
        1 &\  0\\
        1 &\  1\\
        1 &\  1
      \end{array} \right],
    \\
    \tX \in \REAL^{3}_+ \otimes \cC_2 \otimes \REAL^{2}_+ 
    \quad\text{where}\quad 
    &\cC_2 = \operatorname{cone}( W^{(2)} )
    \quad\text{and}\quad 
    W^{(2)} = \left[ \begin{array}{cc}
        1 &\  0\\
        1 &\  0\\
        1 &\  1
      \end{array} \right],
    \\
    \tX \in \REAL^{3}_+ \otimes \REAL^{3}_+ \otimes \cC_3 
    \quad\text{where}\quad 
    &\cC_3 = \operatorname{cone}( W^{(3)} )
    \quad\text{and}\quad 
    W^{(3)} = \left[ \begin{array}{cc}
        1 &\  0\\
        0 &\  1
      \end{array} \right].
  \end{align*}
  
  Recall that in a linear system $AX = B$, if $A$ has full column rank
  then there exists a unique solution $X$.
  Consequently, by taking unfoldings, one finds that if $\tX = \tG
  \times_1 W^{(1)} \times_2 W^{(2)} \times_3 W^{(3)} $ and each
  $W^{(i)}$ has full column rank, then there is a unique solution for
  $\tG$. Note in our example, each $W^{(i)}$ is full column rank.
  Therefore, there is a unique core $\tG$ with the loading matrices
  $W^{(i)}$.
  One can show
  $\matricization_1(\tG) =
  \left[ 
    \begin{array}{cc|cc}
      1 &\ 0 &\ 0 &\ 1\\
      0 &\ 1 &\ 1 &\ -1
    \end{array} 
  \right]$. 
  Since $\tG$ is not nonnegative, by Proposition \ref{prop:ntdCones}
  conclude that $\tX \not\in \cC_1 \otimes \cC_2 \otimes \cC_3$.
  However, if instead the cones corresponding to
  $\bar{W}^{(1)} = \left[ 
    \begin{array}{cc}
      1 &\ 0\\
      0 &\ 1\\
      0 &\ 1
    \end{array} \right]$, 
  $\bar{W}^{(2)} = \left[ 
    \begin{array}{cc}
      1 &\ 0\\
      1 &\ 0\\
      0 &\ 1
    \end{array} 
  \right]$, 
  and $\bar{W}^{(3)} = \left[ \begin{array}{cc}
      1 &\ 0\\
      0 &\ 1
    \end{array} 
  \right]$ 
  were chosen, then $\tX \in \bar{C}^{(1)} \otimes \bar{C}^{(2)}
  \otimes \bar{C}^{(3)}$.
\end{example}

Example \ref{ex:not_in_cones} shows that we cannot take the
intersection of cones to produce a ``minimal'' cone.
Therefore, we make the following mode-wise definition:

\begin{definition}
  Given a nonnegative tensor $\tX \in \REAL^{N_1}_+ \otimes
  \REAL^{N_2}_+ \otimes \REAL^{N_3}_+$, a \emph{minimal $1$-mode
  nonnegative cone}, denoted by $\cC_{1}^{\min}$, is a cone such that
  $\tX\in \cC^{\min}_1\otimes\REAL^{N_2}_+\otimes\REAL^{N_3}_+$ and if
  $\tX\in \cC_1\otimes\REAL^{N_2}_+\otimes\REAL^{N_3}_+$ for some cone
  $\cC_1$, then $\mathcal{O}(\cC^{\min}_1)\leq\mathcal{O}(\cC_1)$
  (we recall that $\mathcal{O}(\cC)$ is the minimum number of vectors
  generating $\cC$).
  We define the $2$-mode and $3$-mode minimal cones analogously.
\end{definition}
Unlike minimal subspaces the minimal cones are generally not unique,
and they are defined mode wise because different mode cones are not
necessarily interchangeable.
We define the \emph{minimal nonnegative multirank} of a nonnegative
tensor as the minimum number of extreme rays of minimal nonnegative
cones along each axis.
\begin{definition}\label{def:nnmultirank2}
  The \emph{$i$-th minimal nonnegative multilinear rank} or
  \emph{$i$-th minimal nonnegative multirank} of a tensor $\tX$,
  denoted $\multirank_{+,i}(\tX)$ is defined as
  $\mathcal{O}(\cC^{\min}_i)$.
  The \emph{minimal nonnegative multilinear rank} of $\tX$ is the
  triple of orders:
  \begin{align*}
    \multirank_{+}(\tX) = (\multirank_{+,1}(\tX), \multirank_{+,2}(\tX), \multirank_{+,3}(\tX)).
  \end{align*}
\end{definition}
As before, we note that the $i$'th minimal nonnegative multilinear
rank does not depend on the $j$'th tensor coordinate for $j \neq i$.
Concretely, we can compute the first minimal nonnegative multirank as
\begin{align*}
  \multirank_{+,1}(\tX) =  
  \min\big\{ \mathcal{O}(\cC_1)~\vert~\tX\in \cC_1 \otimes \REAL_+^{N_2}\otimes \REAL_+^{N_3}, \cC_1 = \operatorname{cone}(W^{(1)}) \subset \REAL_+^{N_1} \big\}
\end{align*}
and, similarly, for $\multirank_{+,2}(\tX)$ and
$\multirank_{+,3}(\tX)$.
As with the real case, it follows directly that for any nTD
$\rank_+(\mF^{(i)}) \geq \multirank_{+,i}(\tX)$. Additionally, we have
an analogous nonnegative statement to Proposition
\ref{prop:unfold_multirank}:

\begin{proposition}
\label{prop:nn:rankunfold_multirank}
  For any nonnegative tensor $\tX$, $\multirank_{+,i}(\tX) =
  \rank_+(\matricization_i(\tX))$.
\end{proposition}
\begin{proof}
  Without loss of generality we prove this for $i=1$ through proving
  the inequality in both directions.
  Let $\multirank_{+,1}(\tX)=k$, then there exists a nonnegative cone
  $\cC^{(1)}$ with $k$ extreme rays such that $\tX \in \cC^{(1)} \otimes
  \REAL_+^{N_2} \otimes \REAL_+^{N_3}$.
  Thus $\tX$ admits a decomposition of the form $\tX = \tD_{identity}
  \times_1 A^{(1)} \times_2 A^{(2)} \times_3 A^{(3)}$ where the
  columns of $A^{(i)}$ are contained by their respective cones,
  $\cC^{(1)}, \REAL^{N_2}_+, \REAL^{N_3}_+$.
  Assemble the extreme rays of $\cC^{(1)}$ into a matrix $W^{(1)} \in
  \REAL^{N_1 \times k}$ so that $\cC^{(1)} =
  \operatorname{cone}(W^{(1)})$.
  Then $A^{(1)} = W^{(1)} H^{(1)}$ for some $H^{(1)} \geq 0$, and with
  substitution we have
  \begin{align*}
    \tX = \tD_{identity} \times_1 W^{(1)}H^{(1)} \times_2 A^{(2)} \times_3 A^{(3)}.
  \end{align*}
  Through distributing and applying unfoldings we have
  \begin{align*}
    \matricization_{1}(\tX) = 
    W^{(1)}\matricization_{1}(\tD_{identity} \times_1 H^{(1)} \times_2 A^{(2)} \times_3 A^{(3)}),
  \end{align*}
  which proves $\rank_+(\matricization_1(\tX)) \leq \rank_+(W^{(1)}) \leq k = 
  \multirank_{+,1}(\tX)$.
  
  Let $\rank_+(\matricization_1(\tX)) = k$.
  Since $\tX$ is nonnegative, $\matricization_1(\tX)$ admits a
  nonnegative decomposition as $\matricization_1(\tX) = W^{(1)}
  H^{(1)}$.
  Since each column of $H^{(1)}$ is nonnegative and is associated with
  a fiber of the tensor, we write the decomposition $\tX =
  \sum_{j=1}^{N_2} \sum_{k=1}^{N_3} W^{(1)}H_{:,j,k} \otimes e_j
  \otimes e_k$.
  This demonstrates that $\tX \in \operatorname{cone}(W^{(1)}) \otimes
  \REAL^{N_2}_+ \otimes \REAL^{N_3}_+$, so $\multirank_{+,i}(\tX) \leq
  \rank_+(\matricization_i(\tX))$.
\end{proof}

Propositions \ref{prop:unfold_multirank} and
\ref{prop:nn:rankunfold_multirank} highlight a key difference between
the real and nonnegative TD.  In the real case, one had that the
minimal subspace was obtained via the unfolding.  In the nonnegative
case, the unfolding does not result in a minimal cone.  From the
definition of the unfolding, one has
\begin{align*}
  \tX 
  = \sum_{j=1}^{N_2}\sum_{k=1}^{N_3}\matricization_1(\tX)_{:,m(j,k)}\otimes\ve_j\otimes\ve_k,
\end{align*}
where $m(j,k)=j+N_{2}(k-1)$, and $\ve_j$ and $\ve_k$ are the $j$-th
and the $k$-th vector of the canonical basis of $\REAL_+^{N_2}$ and
$\REAL_+^{N_3}$, respectively.  Hence,
\begin{align*}
  \tX\in\cone(\matricization_1(\tX))\otimes\REAL_+^{N_2}\otimes\REAL_+^{N_3}.
\end{align*}
However, it may be the case that
$\mathcal{O}(\cone(\matricization_1(\tX))) > \multirank_{1,+}(\tX)$.
Indeed from Proposition \ref{prop:nn:rankunfold_multirank},
$\multirank_{+,i}(\tX) = \rank_+(\matricization_i(\tX))$ and in
general $\mathcal{O}(\cone(\mA)) > \rank_+(\mA)$ for many nonnegative
matrices $\mA$ since the nonnegative rank is equal to the order of the
minimal cone that contains the data.

Just as in the real case, we are interested when the nonnegative TD
has no degeneracy in the loading matrices $\mF^{(i)}$.
When a tensor is simultaneously contained in the tensor product of
minimal nonnegative cones, we call the corresponding nonnegative TD a
\emph{minimal nnTD}.
\begin{definition}
  \label{def:min_nTD}
  An nnTD: $\tX = \tG \times_1 \mF^{(1)} \times_2 \mF^{(2)} \times_3
  \mF^{(3)}$ of a tensor $\tX$ is a \emph{minimal nnTD} whenever the
  core dimensions are equal to the minimal nonnegative multiranks,
  i.e., when $\mF^{(i)}\in\REAL^{N_i\times\multirank_{+,i}(\tX)}$ and
  $\tG\in\REAL^{\multirank_{+,1}(\tX)\times\multirank_{+,2}(\tX)\times\multirank_{+,3}(\tX)}_+$.
\end{definition}
We note that since
$\mF^{(i)}\in\REAL^{N_i\times\multirank_{+,i}(\tX)}$, and
$\rank_+(\mF^{(i)}) \geq \multirank_{+,i}(\tX)$, it follows that
minimal nnTDs satisfy $\rank_+(\mF^{(i)}) =
\multirank_{+,i}(\tX)$. While Example \ref{ex:not_in_cones} showed
that one cannot take intersections to achieve minimal cones, the next
simple result connects the equivalence of simultaneous minimal cones
and a minimal nnTD:

\begin{proposition}
  \label{prop:min_nTD_iff_intersection}
  A nonnegative tensor $\tX$ has a minimal nnTD if and only if there
  exists minimal cones $\cC_i^{\min}$ $i=1,2,3$ for $\tX$ such that $\tX
  \in \cC_1^{\min} \otimes \cC_2^{\min} \otimes \cC_3^{\min}$.
\end{proposition}

\begin{proof}
  Suppose that $\tX$ has a minimal nTD $\tX = \tG \times_1 \mF^{(1)}
  \times_2 \mF^{(2)} \times_3 \mF^{(3)}$ and let $\cC_i =
  \cone(\mF^{(i)})$.  By Proposition \ref{prop:ntdCones}, $\tX \in
  \cC_1 \otimes \cC_2 \otimes \cC_3$. Since
  $\mF^{(i)}\in\REAL^{N_i\times\multirank_{+,i}(\tX)}$ we have
  $\mathcal{O}(\cC_i) \leq \multirank_{+,i}(\tX)$.  However for any
  matrix $\mA$, one has that $\mathcal{O}(\cone(\mA)) \geq
  \rank_+(\mA)$.  Thus $\mathcal{O}(\cC_i) \geq \rank_+(\mF^{(i)}) =
  \multirank_{+,i}(\tX)$.  Combining these two inequalities, we see
  that $\mathcal{O}(\cC_i) = \multirank_{+,i}(\tX)$, so that $\cC_i$
  are minimal order.  Since $\tX \in \cC_1 \otimes \cC_2 \otimes
  \cC_3$, we clearly have $\tX \in \cC_1 \otimes \REAL_+^{N_2} \otimes
  \REAL_+^{N_3}$ so that $\cC_1$ is a minimal cone; likewise for
  $\cC_2$ and $\cC_3$.

  Conversely suppose that $\tX$ has minimal cones $\cC_i^{\min}$
  (namely $\mathcal{O}(\cC_i^{\min}) = \multirank_{+,i}(\tX)$) such
  that $\tX \in \cC_1^{\min} \otimes \cC_2^{\min} \otimes
  \cC_3^{\min}$.  Let $\mF^{(i)}$ be the matrix whose columns are the
  extreme rays of $\cC_i^{\min}$.  Then
  $\mF^{(i)}\in\REAL^{N_i\times\multirank_{+,i}(\tX)}$ and
  $\cone(\mF^{(i)}) = \cC_i^{\min}$.  Since
  \[
  \tX \in \cC_1^{\min} \otimes \cC_2^{\min} \otimes \cC_3^{\min} =  \cone(\mF^{(1)}) \otimes \cone(\mF^{(2)}) \otimes \cone(\mF^{(3)}),
  \]
  by Proposition \ref{prop:ntdCones} there exists a $\tG$ such that
  $\tX = \tG \times_1 \mF^{(1)} \times_2 \mF^{(2)} \times_3
  \mF^{(3)}$.  By construction, this nnTD is a minimal nnTD.
\end{proof}

We further remark from the proof of Proposition
\ref{prop:min_nTD_iff_intersection} that the minimal cones associated
with a minimal nnTD are found by considering the extreme rays of the
cone.

\subsection{Nonnegative Rank Preservation to Tucker Core}

Analogous to the real case it is natural to ask if a minimal nnTD
always exists, or under what conditions does an nnTD exist?
For instance, if $\multirank_{+}(\tX) = (r_1,r_2,r_3)$, then does
there exist nonnegative cones $\cC^{(i)}$ with number of extreme rays
equal to $r_i$ such that $\tX \in \cC^{(1)} \otimes \cC^{(2)} \otimes
\cC^{(3)}$?
Example~\ref{ex:cones_dne}  demonstrates a tensor can fail to have a minimal nnTD:

\begin{example}\label{ex:cones_dne}
  Consider the $4 \times 4 \times 3$ nonnegative tensor with the
  unfoldings,
  \begin{align*}
    \matricization_1(\tX) &= \left[
      \begin{array}{cccc|cccc|cccc}
        1 &\  1 &\  1 &\  1 &\  0 &\  0 &\  1 &\  1 &\  1 &\  2 &\  0 &\  1 \\
        0 &\  1 &\  0 &\  1 &\  0 &\  0 &\  1 &\  1 &\  1 &\  1 &\  0 &\  0 \\
        0 &\  1 &\  0 &\  1 &\  0 &\  1 &\  0 &\  1 &\  0 &\  1 &\  0 &\  1 \\
        1 &\  1 &\  1 &\  1 &\  0 &\  1 &\  0 &\  1 &\  0 &\  2 &\  0 &\  2
      \end{array}
    \right]\;,\\
    \matricization_2(\tX) &= \left[
      \begin{array}{cccc|cccc|cccc}
        1 &\  0 &\  0 &\  1 &\  0 &\  0 &\  0 &\  0 &\  1 &\  1 &\  0 &\  0 \\
        1 &\  1 &\  1 &\  1 &\  0 &\  0 &\  1 &\  1 &\  2 &\  1 &\  1 &\  2 \\
        1 &\  0 &\  0 &\  1 &\  1 &\  1 &\  0 &\  0 &\  0 &\  0 &\  0 &\  0 \\
        1 &\  1 &\  1 &\  1 &\  1 &\  1 &\  1 &\  1 &\  1 &\  0 &\  1 &\  2
      \end{array}
    \right]\;,\\
    \matricization_3(\tX) &= \left[
      \begin{array}{cccc|cccc|cccc|cccc}
        1 &\  0 &\  0 &\  1 &\  1 &\  1 &\  1 &\  1 &\  1 &\  0 &\  0 &\  1 &\  1 &\  1 &\  1 &\  1 \\
        0 &\  0 &\  0 &\  0 &\  0 &\  0 &\  1 &\  1 &\  1 &\  1 &\  0 &\  0 &\  1 &\  1 &\  1 &\  1 \\
        1 &\  1 &\  0 &\  0 &\  2 &\  1 &\  1 &\  2 &\  0 &\  0 &\  0 &\  0 &\  1 &\  0 &\  1 &\  2 
      \end{array}
    \right]\;.
  \end{align*}
  Suppose there exists a minimal nTD
  \begin{align*}
    \tX = \tG \times_1 \mF^{(1)} \times_2 \mF^{(2)} \times_3 \mF^{(3)}
  \end{align*}
  with
  $\tG\in\REAL_+^{\multirank_{+,1}(\tX)\times\multirank_{+,2}(\tX)\times\multirank_{+,3}(\tX)}$.
  From the decomposition  
  \begin{align*}
    \matricization_2(\tX) = \left[
      \begin{array}{ccc}
        0 &\  0 &\  1 \\
        0 &\  1 &\  1 \\
        1 &\  0 &\  0 \\
        1 &\  1 &\  0
      \end{array}
    \right] \left[
      \begin{array}{cccc|cccc|cccc}
        1 &\  0 &\  0 &\  1 &\  1 &\  1 &\  0 &\  0 &\  0 &\  0 &\  0 &\  0 \\
        0 &\  1 &\  1 &\  0 &\  0 &\  0 &\  1 &\  1 &\  1 &\  0 &\  1 &\  2 \\
        1 &\  0 &\  0 &\  1 &\  0 &\  0 &\  0 &\  0 &\  1 &\  1 &\  0 &\  0
      \end{array}
    \right]\;
  \end{align*}
  it can be verified that $\multirank_{+,2}(\tX) = 3$,
  and therefore $\mF^{(2)} \in \REAL^{4 \times 3}_+$.  
  From the decomposition
  \begin{align*}
    \matricization_1(\tX) = \left[
      \begin{array}{ccc}
        0 &\  1 &\  1 \\
        0 &\  0 &\  1 \\
        1 &\  0 &\  0 \\
        1 &\  1 &\  0
      \end{array}
    \right] \left[
      \begin{array}{cccc|cccc|cccc}
        0 &\  1 &\  0 &\  1 &\  0 &\  1 &\  0 &\  1 &\  0 &\  1 &\  0 &\  1 \\
        1 &\  0 &\  1 &\  0 &\  0 &\  0 &\  0 &\  0 &\  0 &\  1 &\  0 &\  1 \\
        0 &\  1 &\  0 &\  1 &\  0 &\  0 &\  1 &\  1 &\  1 &\  1 &\  0 &\  0 
      \end{array}
    \right]
  \end{align*}
  it can be verified that $\multirank_{+,1}(\tX) = 3$, and therefore
  $\mF^{(1)} \in \REAL^{4 \times 3}_+$.
  This decomposition, and the corresponding tensor decomposition
  $\tX=\tH\times_1 W$, where
  \begin{align*}
    W = \left[ \begin{array}{ccc}
        0 &\  1 &\  1 \\
        0 &\  0 &\  1 \\
        1 &\  0 &\  0 \\
        1 &\  1 &\  0
      \end{array} \right], \matricization_1(\tH) = \left[
      \begin{array}{cccc|cccc|cccc}
        0 &\  1 &\  0 &\  1 &\  0 &\  1 &\  0 &\  1 &\  0 &\  1 &\  0 &\  1 \\
        1 &\  0 &\  1 &\  0 &\  0 &\  0 &\  0 &\  0 &\  0 &\  1 &\  0 &\  1 \\
        0 &\  1 &\  0 &\  1 &\  0 &\  0 &\  1 &\  1 &\  1 &\  1 &\  0 &\  0 
      \end{array}
    \right],
  \end{align*}
  is unique (\cite{Gillis:2012} Theorem 6), which implies that with proper permutation and scaling $W = \mF^{(1)}$, and $\tH = \tG
  \times_2 \mF^{(2)} \times_3 \mF^{(3)}$.
  Note that the $\multirank_{+,2}(\tH) = 4$ since the second unfolding
  \begin{align*}
    \matricization_2(\tH) = \left[
      \begin{array}{ccc|ccc|ccc}
        0 &\  1 &\  0 &\  0 &\  0 &\  0 &\  0 &\  0 &\  1 \\
        1 &\  0 &\  1 &\  1 &\  0 &\  0 &\  1 &\  1 &\  1 \\
        0 &\  1 &\  0 &\  0 &\  0 &\  1 &\  0 &\  0 &\  0 \\
        1 &\  0 &\  1 &\  1 &\  0 &\  1 &\  1 &\  1 &\  0 
      \end{array}
    \right],
  \end{align*}
  contains Example~\ref{ex:mrank_mrank+} as a submatrix. 
  However, by Proposition~\ref{prop:nn:rankunfold_multirank}
  \begin{align*}
    \multirank_{+,2}(\tH) &\ = \multirank_{+,2}(\tG \times_2 \mF^{(2)} \times_3 \mF^{(3)}) \\
    &\ = \rank_+(\matricization_2(\tG \times_2 \mF^{(2)} \times_3 \mF^{(3)})\\
    &\ = \rank_+(\mF^{(2)} \matricization_2(\tG \times_3 \mF^{(3)}))\leq 3.
  \end{align*}
  This is a contradiction, so the supposition that there exists a
  minimal nnTD is false.
\end{example}

A further question is: if the nnTD: $\tX=\tG\times_1 \mF^{(1)}\times_2
\mF^{(2)}\times_3 \mF^{(3)}$ does exist, is the nonnegative rank of $\tX$
preserved to the nnTD, $\tG$, that is, is $\rank_+(\tX) =
\rank_+(\tG)?$
Example~\ref{ex:equal} demonstrates that even when the minimal nnTD does exist,
the nonnegative rank of the tensor is not necessarily preserved to the
core.

\begin{example}\label{ex:equal}
  Let $\tX = \tD_{identity} \times_1 A^{(1)}_{\tX} \times_2
  A^{(2)}_{\tX} \times_3 A^{(3)}_{\tX}$, where $\tD_{identity}$ is the
  diagonal identity tensor and
  \begin{align*}
    A^{(1)}_{\tX} = \begin{bmatrix} 
      1 &\  1 &\  0 &\  0 \\ 
      1 &\  0 &\  1 &\  0 \\ 
      0 &\  1 &\  0 &\  1 \\ 
      0 &\  0 &\  1 &\  1
    \end{bmatrix},
    \quad
    A^{(2)}_{\tX} = \begin{bmatrix} 
      1 &\  1 &\  1 &\  1 \\ 
      0 &\  1 &\  1 &\  1 \\ 
      0 &\  0 &\  1 &\  1 \\ 
      0 &\  0 &\  0 &\  1  
    \end{bmatrix}, 
    \quad
    A^{(3)}_{\tX} = \begin{bmatrix} 
      1 &\  1 &\  1 &\  1 \\ 
      0 &\  1 &\  1 &\  1 \\ 
      0 &\  0 &\  1 &\  1 \\ 
      0 &\  0 &\  0 &\  1  
    \end{bmatrix}.
  \end{align*}
  Kruskal's theorem~\cite{Kruskal:1977,Domanov-DeLathauwer:2013}
  proves that $\rank_+(\tX)=4$, and the nnCPD of the tensor is unique.
  Using Proposition \ref{prop:nn:rankunfold_multirank}, one can show
  that $\multirank_+(\tX) = (3,4,4)$. For example, from the first
  unfolding of $\tX$, we have
  \begin{align*}
    \matricization_1(\tX) = \left[
      \begin{array}{cccc|cccc|cccc|cccc}
        2 &\  1 &\  0 &\  0 &\  1 &\  1 &\  0 &\  0 &\  0 &\  0 &\  0 &\  0 &\  0 &\  0 &\  0 &\  0 \\
        2 &\  1 &\  1 &\  0 &\  1 &\  1 &\  1 &\  0 &\  1 &\  1 &\  1 &\  0 &\  0 &\  0 &\  0 &\  0 \\
        2 &\  2 &\  1 &\  1 &\  2 &\  2 &\  1 &\  1 &\  1 &\  1 &\  1 &\  1 &\  1 &\  1 &\  1 &\  1 \\
        2 &\  2 &\  2 &\  1 &\  2 &\  2 &\  2 &\  1 &\  2 &\  2 &\  2 &\  1 &\  1 &\  1 &\  1 &\  1
      \end{array}
    \right].
  \end{align*}
  The rank of $\matricization_1(\tX)$ is $3$, and this matrix admits a
  nonnegative decomposition,
  \begin{align*}
    \matricization_1(\tX) = \begin{bmatrix} 
      0 &\  0 &\  1 \\ 
      0 &\  1 &\  1 \\ 
      1 &\  0 &\  0 \\ 
      1 &\  1 &\  0
    \end{bmatrix}
    \left[
      \begin{array}{cccc|cccc|cccc|cccc}
        2 &\  2 &\  1 &\  1 &\  2 &\  2 &\  1 &\  1 &\  1 &\  1 &\  1 &\  1 &\  1 &\  1 &\  1 &\  1 \\
        0 &\  0 &\  1 &\  0 &\  0 &\  0 &\  1 &\  0 &\  1 &\  1 &\  1 &\  0 &\  0 &\  0 &\  0 &\  0 \\
        2 &\  1 &\  0 &\  0 &\  1 &\  1 &\  0 &\  0 &\  0 &\  0 &\  0 &\  0 &\  0 &\  0 &\  0 &\  0
      \end{array}
    \right].
  \end{align*}
  Using these decompositions, one can show that $\tX$ admits a minimal
  nnTD of the form $\tX = \tG \times_1 \mF^{(1)} \times_2 I \times_3 I$
  where
  \begin{align*}
    \matricization_1(\tG) =  \left[    
      \begin{array}{cccc|cccc|cccc|cccc}
        2 &\  2 &\  1 &\  1 &\  2 &\  2 &\  1 &\  1 &\  1 &\  1 &\  1 &\  1 &\  1 &\  1 &\  1 &\  1 \\
        0 &\  0 &\  1 &\  0 &\  0 &\  0 &\  1 &\  0 &\  1 &\  1 &\  1 &\  0 &\  0 &\  0 &\  0 &\  0 \\
        2 &\  1 &\  0 &\  0 &\  1 &\  1 &\  0 &\  0 &\  0 &\  0 &\  0 &\  0 &\  0 &\  0 &\  0 &\  0
      \end{array}
    \right]
  \end{align*}
  and 
  \begin{align*}
    \mF^{(1)} = \begin{bmatrix} 
      0 &\  0 &\  1 \\ 
      0 &\  1 &\  1 \\ 
      1 &\  0 &\  0 \\ 
      1 &\  1 &\  0
    \end{bmatrix}.
  \end{align*}
  Now, suppose to the contrary that $\rank_+(\tX) = \rank_+(\tG)$.
  Then, let $\tG = \tD_{\tG} \times_1 \mA^{(1)}_{\tG} \times_2
  A^{(2)}_{\tG} \times_3 A^{(3)}_{\tG}$ be an nCPD of $\tG$.
  Since the nCPD of $\tX$ is unique, we have up to permutation and
  nonnegative scaling that
  \begin{align}\label{ex:eq:matrix_decomposition}
    A^{(1)}_{\tX} = \mF^{(1)} \mA^{(1)}_{\tG}.
  \end{align}
  From Example~\ref{ex:mrank_mrank+} we know that $\rank_+(\mA^{(1)}_{\tX}) = 4$.  But then
  \begin{align*}
    \rank_+(A^{(1)}_{\tX}) > 3 = \rank_+(\mF^{(1)}) \geq \rank_+(\mF^{(1)} \mA^{(1)}_{\tG}),
  \end{align*}
  which is a contradiction. 
  Therefore, $\rank_+(\tX) \neq \rank_+(\tG)$.
\end{example}

Example \ref{ex:equal} highlights a key difference between the real
and nonnegative minimal Tucker decompositions. By Theorem
\ref{def:cpd-td}, if the CPD is unique then the column space of the
CPD loading matrices will recover the minimal subspaces.  However when
the nonnegative CPD is unique, the loading matrices can still fail to
capture the minimal cone ($\cC^{(1)}_{\min}$ and $A^{(1)}_{\tX}$ in
previous example). In particular, the nnCPD cannot be a minimal Tucker
decomposition in this case. This causes problems with preservation of
the rank to the core of the tensor.  It turns out that this issue in
Example \ref{ex:equal} is always hold. Namely, when $\tX$ has a unique
nnCPD and a Tucker has a factor with nonnegative rank smaller then the
loading matrix factor in the nnCPD, nonnegative rank cannot be
preserved:

\begin{theorem}\label{thm: rank not equal}
  Let $\tX$ be a nonnegative tensor with unique nCPD $\tX =
  \tD_{identity} \times_1 \mA^{(1)}_{\tX} \times_2 \mA^{(2)}_{\tX}
  \times_3 \mA^{(3)}_{\tX}$.  Suppose that $\tX$ has a nTD $\tX = \tG
  \times_1 \mF^{(1)} \times_2 \mF^{(2)} \times_3 \mF^{(3)}$ where
  $\rank_+(\mA^{(i)}) > \rank_+(\mF^{(i)})$ for some $i=1,2,3$.  Then
  $\rank_+(\tX) \neq \rank_+(\tG)$.
\end{theorem}

\begin{proof}
  Without loss of generality, let $i=1$, i.e. $\rank_+(\mA^{(1)}) >
  \rank_+(\mF^{(1)})$. Suppose to the contrary that $\rank_+(\tX) =
  \rank_+(\tG)$.  Let $\tG = \tD_{\tG} \times_1 \mA^{(1)}_{\tG}
  \times_2 \mA^{(2)}_{\tG} \times_3 \mA^{(3)}_{\tG}$ be an nCPD of
  $\tG$.  Then both $\tX = \tD_{identity} \times_1 \mA^{(1)}_{\tX}
  \times_2 \mA^{(2)}_{\tX} \times_3 \mA^{(3)}_{\tX}$ and
  \[
  \tX = \tG \times_1 \mF^{(1)} \times_2 \mF^{(2)} \times_3 \mF^{(3)} = \tD_{\tG} \times_1 (\mF^{(1)} \mA^{(1)}_{\tG}) \times_2 (\mF^{(1)} \mA^{(2)}_{\tG}) \times_3 (\mF^{(1)} \mA^{(3)}_{\tG})
  \]
  are rank decompositions of $\tX$. Since the nCPD of $\tX$ is unique,
  up to permutation and nonnegative scaling one has
  \[
    \mA^{(1)}_{\tX} = \mF^{(1)} \mA^{(1)}_{\tG}.
  \]
  Thus, $\rank_+(\mA^{(1)}_{\tX}) = \rank_+(\mF^{(1)} \mA^{(1)}_{\tG})$. However, by assumption 
  \[
  \rank_+(\mA^{(1)}_{\tX}) > \rank_+(\mF^{(1)}) \geq \rank_+(\mF^{(1)} \mA^{(1)}_{\tG}),
  \]
  a contradiction. 
\end{proof}

Theorem \ref{thm: rank not equal} gives condition on when the rank is
\textit{not} preserved based on the CPD.  We believe that for a large
class of nonnegative tensors where compression is achieved, this
implies that the rank is not preserved. However, deriving precises
probabilistic statements is challenging due to the non-stochastic
relationship between loading matrices and random tensors.

Examples \ref{ex:cones_dne} and \ref{ex:equal} demonstrate the
subtleties of the nonnegative factorizations compared to the real
valued.  First, the minimal nnTD can fail to exist.
Second, even if it exists, the nonnegative rank of the minimal nnTD
core may not be equal to the nonnegative rank of $\tX$.
The following Theorem provides sufficient conditions for a minimal
nnTD to exist, and for the nonnegative rank of the tensor to be
preserved to the core of the minimal nnTD.  We note that because of
Theorem \ref{thm: rank not equal}, a rank requirement for nnCPD
loading matrices is required. The following is the nonnegative analog
of Theorem \ref{thm:minTDExists}.

\begin{theorem}\label{thm:minnTDExists}
  Suppose a nonnegative tensor $\tX$ has an nCPD:
  $\tX=\tD_{\tX}\times_1 \mA_{\tX}^{(1)}\times_2
  \mA_{\tX}^{(2)}\times_3 \mA_{\tX}^{(3)}$ with
  $\rank_+(\mA_{\tX}^{(i)}) = \multirank_{+,i}(\tX)$ for $1 \leq i
  \leq 3$.  Then a minimal nTD: $\tX = \tG \times_1 \mF^{(1)} \times_2
  \mF^{(2)} \times_3 \mF^{(3)}$ exists such that $\rank_+(\tX) =
  \rank_+(\tG)$.
\end{theorem}
\begin{proof}
  Since $\rank_+(\mA_{\tX}^{(i)}) = \multirank_{+,i}(\tX)$ for $1 \leq
  i \leq 3$, each $\mA_{\tX}^{(i)}$ has a nonnegative decomposition
  $\mA_{\tX}^{(i)} = W^{(i)}H^{(i)}$.
  Substituting into the nCPD and distributing
  \begin{align*}
    \tX 
    & = \tD_{\tX} \times_1 \mA_{\tX}^{(1)} \times_2 \mA_{\tX}^{(2)} \times_3 \mA_{\tX}^{(3)}\\
    &= \tD_{\tX} \times_1 \mW^{(1)}\mH^{(1)} \times_2 \mW^{(2)}\mH^{(2)} \times_3 \mW^{(3)}\mH^{(3)}\\
    &= (\tD_{\tX} \times_1 \mH^{(1)} \times_2 \mH^{(2)} \times_3 \mH^{(3)}) \times_1 \mW^{(1)} \times_2 \mW^{(2)} \times_3 \mW^{(3)}\\
    &= \tG \times_1 \mW^{(1)} \times_2 \mW^{(2)} \times_3 \mW^{(3)}
  \end{align*}
  where $\tG=\tD_{\tX}\times_1 \mH^{(1)}\times_2 \mH^{(2)}\times_3
  \mH^{(3)}$.
  The core $\tG$ is a nonnegative tensor with shape equal to the
  nonnegative minimal  multiranks of $\tX$, so $\tX = \tG \times_1 \mW^{(1)}
  \times_2 \mW^{(2)} \times_3 \mW^{(3)}$ is a minimal nTD.
  To prove $\rank_+(\tX) = \rank_+(\tG)$, it once again suffices to show that $\rank_+(\tX) \geq \rank_+(\tG)$.
  However from the constructed decomposition $\tG=\tD_{\tX}\times_1 \mH^{(1)}
  \times_2 \mH^{(2)}\times_3 \mH^{(3)}$ we know
  $\rank_+(\tG)\leq\rank_+(\tX)$.
\end{proof}

Theorem \ref{thm:minnTDExists} demonstrates that a minimal nnTD exists
that will preserve the rank.  However contrary to Theorem
\ref{thm:minTDExists}, it does not state that every minimal nnTD will
preserve the rank to the core. This is yet another fundamental
challenge one must surmount in the nonnegative case - not every
minimal nnTD will necessarily preserve the rank.  The following
example illustrates this issue:

\begin{example}
  \label{ex:not_every_min_preserves_rank}
  Let $\tX \in \REAL_+^{3,3,3}$ be the tensor given by 
  \[
  \tX = 
  \left[    
    \begin{array}{ccc|ccc|ccc}
      2 &\ 8 &\ 3 &\ 1 &\ 5 &\ 2 &\ 2 &\ 8 &\ 3\\
      4 &\ 15 &\ 5 &\ 2 &\ 8 &\ 3 &\ 4 &\ 15 &\ 5\\
      2 &\ 6 &\ 2 &\ 1 &\ 3 &\ 1 &\ 2 &\ 6 &\ 2
    \end{array}
    \right].
  \]
  Then $\tX$ has an nnCPD given by
  \begin{align*}
    A^{(1)}_{\tX} = \begin{bmatrix} 
      2 & 1 & 1\\
      1 & 1 & 0\\
      2 & 1 & 1
    \end{bmatrix},
    \quad
    A^{(2)}_{\tX} = \begin{bmatrix} 
      1 & 1 & 0\\
      2 & 1 & 1\\
      1 & 0 & 0
    \end{bmatrix}, 
    \quad
    A^{(3)}_{\tX} = \begin{bmatrix} 
      1 & 0 & 0\\
      3 & 2 & 1\\
      1 & 1 & 0
    \end{bmatrix}.
  \end{align*}
  One can check that $\rank_+(A^{(1)}_{\tX}) = 2 =
  \multirank_{+,1}(\tX)$ and $\rank_+(A^{(2)}_{\tX}) =
  \rank_+(A^{(3)}_{\tX}) = 3 = \multirank_{+,2}(\tX) =
  \multirank_{+,3}(\tX)$. Thus, $\tX$ satisfies the hypothesis of
  Theorem \ref{thm:minnTDExists}. We will now show that there are two
  minimal nTDs
  \[
  \tX = \tG_1 \times_1 \mF_1^{(1)} \times_2 \mF_1^{(2)} \times_3 \mF_1^{(3)} = \tG_2 \times_1 \mF_2^{(1)} \times_2 \mF_2^{(2)} \times_3 \mF_2^{(3)}
  \]
  with $\rank(\tX) = \rank_+(\tG_1) < \rank_+(\tG_2)$.  Therefore, not
  every minimal nnTD can preserve the rank to the Tucker core. Indeed,
  one can check
  \begin{align*}
    \tG_1 = 
    \left[    
      \begin{array}{ccc|ccc}
        1 & 1 & 1 & 1 & 1 & 1\\
        1 & 2 & 2 & 1 & 1 & 1\\
        1 & 1 & 1 & 1 & 2 & 1
      \end{array}
      \right],
    \quad
    \mF_1^{(1)} = \begin{bmatrix} 
      1 & 1\\
      1 & 0\\
      1 & 1
    \end{bmatrix},
    \quad
    \mF_1^{(2)} = \begin{bmatrix} 
      0 & 1 & 0\\
      0 & 1 & 1\\
      1 & 0 & 0
    \end{bmatrix}, 
    \quad
    \mF_1^{(3)} = \begin{bmatrix} 
      1 & 0 & 0\\
      1 & 1 & 1\\
      0 & 0 & 1
    \end{bmatrix}
  \end{align*}
  and 
  \begin{align*}
    \tG_2 = 
    \left[    
      \begin{array}{ccc|ccc}
        0 & 2 & 1 & 1 & 3 & 1\\
        0 & 1 & 1 & 2 & 7 & 2\\
        0 & 0 & 0 & 1 & 3 & 1
      \end{array}
      \right],
    \quad
    \mF_2^{(1)} = \begin{bmatrix} 
      1 & 2\\
      1 & 1\\
      1 & 2
    \end{bmatrix},
    \quad
    \mF_2^{(2)} = \begin{bmatrix} 
      1 & 0 & 0\\
      0 & 1 & 0\\
      0 & 0 & 1
    \end{bmatrix}, 
    \quad
    \mF_2^{(3)} = \begin{bmatrix} 
      1 & 0 & 0\\
      0 & 1 & 0\\
      0 & 0 & 1
    \end{bmatrix}
  \end{align*}
  result in minimal nTDs of $\tX$. From the decomposition
  \[
  \tG_1 = 
  \left[    
    \begin{array}{ccc|ccc}
      1 & 1 & 1 & 1 & 1 & 1\\
      1 & 1 & 1 & 1 & 1 & 1\\
      1 & 1 & 1 & 1 & 1 & 1
    \end{array}
    \right]
  +
  \left[    
    \begin{array}{ccc|ccc}
      0 & 0 & 0 & 0 & 0 & 0\\
      0 & 1 & 1 & 0 & 0 & 0\\
      0 & 0 & 0 & 0 & 0 & 0\\
    \end{array}
    \right]
  +
  \left[    
    \begin{array}{ccc|ccc}
      0 & 0 & 0 & 0 & 0 & 0\\
      0 & 0 & 0 & 0 & 0 & 0\\
      0 & 0 & 0 & 0 & 1 & 0\\
    \end{array}
    \right]
  \]
  we see that $\rank_+(\tG_1) \leq 3$.  Since $\rank_+(\tX) = 3$ by
  above, we have that $\rank_+(\tG_1) = 3$. We now show that
  $\rank_+(\tG_2) \neq 3$. Indeed by Kruskal's
  Theorem~\cite{Kruskal:1977}, $\tX$ has a unique nnCPD. By Proposition
  1 of ~\cite{Cohen-Comon-Gillis:2017}, $\rank_+(\tG_1) = \rank_+(\tX) =
  3$ if and only if $A^{(i)}_{\tX} \subset \cone(\mF_2^{(i)})$. Since
  $[1,0,1] \notin \cone(\mF_2^{(1)})$, we see that the rank cannot be
  preserved to the core $\tG_2$.
\end{example}

Cohen~\textit{et al.}~\cite{Cohen-Comon-Gillis:2017} (see
Proposition~1) provide some necessary and sufficient conditions for
the nonnegative rank of a tensor to persist to the core of an nnTD
under some geometric hypothesis.  We remark that their theorem, as
stated, requires uniqueness of the nnCPD along with a full column rank
condition on the factors of the nnCPD. However the full column rank is
not needed, and the uniqueness of the nnCPD is only required for one
direction.  Namely, that if the nnCPD is unique and the rank of the
tensor is preserved to the core, then the (unique) nnCPD factors are
contained inside the cones from the nnTD loading matrices. We made the
equivalent converse statement in Theorem \ref{thm:nCPDmultirank}.  It
too, requires uniqueness of the nnCPD.  However, since one half of our
`if and only if' does not require uniqueness, we have opted to
separate the two conditions.

\subsection{Nonnegative Rank Deficiency in nnCPD Factors}

We now discuss the difficulties associated from nonnegative rank
deficent factors in a nnCPD.  We recap the work above to discuss the
challenges facing a nonnegative analog of CANDELINC - in particular,
the issues surrounding the existence of a min nnTD. Following the
subsection on real valued CANDELINC above, we begin by exploring the
relations between the nonnegative minimal multirank and the
nonnegative ranks of nnCPD factors.  The following theorem relates the
uniqueness of a nnCPD to the minimal nnTDs:

\begin{theorem}\label{thm:nCPDmultirank}
  Let $\tX=\tD\times_1 \mA^{(1)}\times_2 \mA^{(2)}\times_3 \mA^{(3)}$ be an
  nnCPD, then $\multirank_{+,i}(\tX)\leq\rank_+(\mA^{(i)})$. 
  Furthermore, if the nnCPD is unique and there exists a minimal nnTD
  with $\rank_+(\tX) = \rank_+(\tG)$ then $\multirank_{+,i}(\tX) =
  \rank_+(\mA^{(i)})$
\end{theorem}
\begin{proof}
  That $\multirank_{+,i}(\tX)\leq\rank_+(\mA^{(i)})$ follows from the
  fact that an nnCPD is a nnTD.  Now suppose that $\tX$ has a unique
  nCPD, and consider a minimal nTD of the form
  $\tX=\tG\times_1\mF^{(1)}\times_2\mF^{(2)}\times_3\mF^{(3)}$.
  We recall that $\rank_+(\mF^{(i)})=\multirank_{+,i}(\tX)$ by
  definition.
  On its turn, the Tucker core admits an nCPD
  $\tG=\tD_{\tG}\times_1\mB^{(1)}\times_2\mB^{(2)}\times_3\mB^{(3)}$.
  Substituting the nCPD of $\tG$ in the nTD of $\tX$, we obtain the
  alternative nCPD
  $\tX=\tD_{\tG}\times_1\mF^{(1)}\mB^{(1)}\times_2\mF^{(2)}\mB^{(2)}\times_3\mF^{(3)}\mB^{(3)}$.
  Since we assume that the nCPD is unique, with appropriate
  nonnegative scalings and permutations, which are nonnegative
  rank-preserving operations, we obtain that
  $\mA^{(i)}=\mF^{(i)}\mB^{(i)}$, which proves that
  $\rank_+(\mA^{(i)})\leq\rank_+(\mF^{(i)})=\multirank_{+,i}(\tX)$.
\end{proof}
Theorem~\ref{thm:nCPDmultirank} illuminates an algorithmic challenge
of computing a nCPD directly.  Suppose $\tX
=\tD\times_1\mA^{(1)}\times_2\mA^{(2)}\times_3\mA^{(3)}$ has a unique
nCPD with $\rank_+(\tX) = r$ and minimal multiranks $\multirank_i(\tX)
= r_i$. In practice, many shaped tensors have $r_i<r$.  Thus,
$\mA^{(i)}$ is a rank deficient matrix by Theorem
\ref{thm:nCPDmultirank}, which can be an algorithmically challenging
task to overcome.

This motivates the need for nonnegative version of CANDELINC.
Concretely, if
$\tX=\tG\times_1\mF^{(1)}\times_2\mF^{(2)}\times_3\mF^{(3)}$ is a
minimal nnTD such that $\rank_+(\tG) = \rank_+(\tX)$, and
$\tG=\tD_{\tG}\times_1\mA_{\tG}^{(1)}\times_2\mA_{\tG}^{(2)}\times_3\mA_{\tG}^{(3)}$
is the CPD of the Tucker core, then each factor $\mA_{\tG}^{(i)}$ is a
full column rank matrix This suggests that, under some conditions, a
stable way of computing nnCPD is to first compute a minimal nnTD,
e.g., using approximate NMF, then compute nCPD on the nonnegative
Tucker core, and finally substitute the nnCPD of the core-tensor in
the nnTD to obtain the final nnCPD of the original tensor.

Unfortunately, the previous work highlights some of the major
challenges a nonnegative CANDELINC must overcome.  Indeed a min nnTD
need not exist (Example \ref{ex:cones_dne}), and even when it does, it
need not preserve the rank to the core (Example \ref{ex:equal}).
Furthermore, a tensor may have a minimal nnTD which preserves the
rank, but this does not mean all minimal nnTD will preserve the rank
to the core (Example \ref{ex:not_every_min_preserves_rank}).

These issues indicate that minimal Tucker decompositions, while
desirable, are perhaps not feasible in the nonnegative case.  To
overcome this hurdle, we redirect our interest to nonnegative Tucker
decompostions we call \textit{canonical}. In the next subsection, we
define the canonical Tucker decomposition and show that one always
exists which preserves the rank to the core.

\subsection{The Canonical Multirank and Canonical Tucker}

Ultimately, we desire a decomposition that 1) preserves the nCPD rank
to the core and 2) does not require extraction of rank deficient
matrices. While minimal nnTD's are devoid of rank deficiency in the
loading matrices, they may not preserve the rank to the core (if they
exist at all).  This conversation leads us to desire a less strict
type of Tucker decomposition which not only preservers the rank, but
whose loading matrices are also devoid of rank deficiency. We
therefore make the following definition:

\begin{definition}
  \label{def: canonical rank}
  Let $\tX \in \REAL_+^{N_1\times N_2 \times N_3}$ have a unique nnCPD
  given by
  $\tX=\tD\times_1\mA^{(1)}\times_2\mA^{(2)}\times_3\mA^{(3)}$.  The
  \textit{i-th canonical (nonnegative) multirank of $\tX$}, denoted
  $\kappa_i(\tX)$, is the rank of the i-th nnCPD factor
  $\rank_+(\mA^{(i)})$. The \textit{canonical (nonnegative) multirank}
  of $\tX$ is the triple
  \[
  \kappa(\tX) = (\kappa_1(\tX), \kappa_2(\tX), \kappa_3(\tX)) = (\rank_+(\mA^{(1)}), \rank_+(\mA^{(2)}), \rank_+(\mA^{(3)})).
  \]
\end{definition}

Corresponding to the canonical multirank we have a canonical Tucker decomposition. 

\begin{definition}
  \label{def: canonical tucker}
    Consider tensor $\tX\in\REAL_+^{N_1\times N_2\times N_3}$ with unique nnCPD.
  We say that the Tucker decomposition
  $\tX=\tG\times_1\mF^{(1)}\times_2\mF^{(2)}\times_3\mF^{(3)}$ is
  \emph{canonical} if the dimensions of the core tensor $\tG$ are equal
  to the canonical multiranks, i.e.,
  $\tG\in\REAL_+^{\kappa_1(\tX)\times\kappa_2(\tX)\times\kappa_3(\tX)}$
  and $\mF^{(i)}\in\REAL_+^{N_i\times\kappa_i(\tX)}$, $i=1,2,3$.
\end{definition}

By definition $\multirank_{i,+}(\tX) \leq \kappa_i(\tX)$.  Therefore,
a canonical Tucker is a less restrictive shape constraint than a
minimal Tucker.  Unlike minimal Tuckers, one can always find a
canonical Tucker that will preserve the rank:

\begin{theorem}
\label{thm: canonical preserve rank}
  Let $\tX \in \REAL_+^{N_1\times N_2 \times N_3}$ have unique
  nnCPD. Then there exists a canonical nnTD which preserves the rank.
\end{theorem}

\begin{proof}
 Let $\tX \in \REAL_+^{N_1\times N_2 \times N_3}$ have unique nnCPD
 given by
 $\tX=\tD\times_1\mA^{(1)}\times_2\mA^{(2)}\times_3\mA^{(3)}$. For
 each $i=1,2,3$, consider the rank factorizations $\mA^{(i)} =
 \mW^{(i)} \mH^{(i)}$.  Subbing these factorizations into the nnCPD:
 \[
 \begin{array}{rcl}
    \tX & = & \tD\times_1\mA^{(1)}\times_2\mA^{(2)}\times_3\mA^{(3)}\\
    & = & \left(\tD\times_1\mH^{(1)}\times_2\mH^{(2)}\times_3\mH^{(3)} \right) \times_1\mW^{(1)}\times_2\mW^{(2)}\times_3\mW^{(3)}.
 \end{array}
 \]
 Let $\tG :=
 \tD\times_1\mH^{(1)}\times_2\mH^{(2)}\times_3\mH^{(3)}$. By
 construction, $\mW^{(i)} \in \REAL_+^{N_i\times \kappa_i(\tX)}$ so
 that $\tX = \tG \times_1\mW^{(1)}\times_2\mW^{(2)}\times_3\mW^{(3)} $
 is a canonical Tucker. Furthermore, $\rank_+(\tG) \leq \rank_+(\tX)$
 so that $\rank_+(\tG) = \rank_+(\tX)$.
\end{proof}

Theorem \ref{thm: canonical preserve rank} states that there exists a
rank preserving canonical Tucker.  However, Example
\ref{ex:not_every_min_preserves_rank} shows that not every canonical
Tucker can preserve the rank. Unlike its real counterparts, the shape
of a nonnegative Tucker decomposition does not guarantee that all such
factorizations will preserve rank.

In minimal nnTD, the rank of the loading matrix $\mF^{(i)}$ is equal
to $\multirank_{i,}(\tX)$ so that the matrix is not degenerate. Note
that in a canonical nnTD
$\tX=\tG\times_1\mF^{(1)}\times_2\mF^{(2)}\times_3\mF^{(3)}$, one has
that $\multirank_{i,+}(\tX) \leq \rank_+(\mF^{(i)}) \leq
\kappa_i(\tX)$.  If one selects a canonical nnTD which preserves the
rank to the core, then $\rank_+(\mF^{(i)}) = \kappa_i(\tX)$ so that
once again, the matrix is not degenerate. The challenge here lies
instead in the selection of the correct cones that contains the
minimal cones needed to preserve the rank.

    
    
    
    

\begin{figure}[!t]
  \centering
  \includegraphics[width=0.365\textwidth]{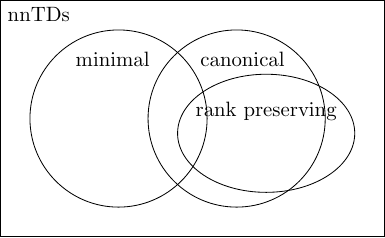}
  \caption{Diagram of nonnegative Tucker decompositions for tensors with a unique nnCPD.}
  \label{fig:nTDschematic}
\end{figure}

\begin{proposition}
  \label{prop: canonical}
  Let $\tX \in \REAL_+^{N_1\times N_2 \times N_3}$ have unique nnCPD
  given by
  $\tX=\tD\times_1\mA^{(1)}\times_2\mA^{(2)}\times_3\mA^{(3)}$.
  Suppose $\tX=\tG\times_1\mF^{(1)}\times_2\mF^{(2)}\times_3\mF^{(3)}$
  is a canonical Tucker such that $R = \rank_+(\tX) = \rank_+(\tG)$.
  Then $\rank_+(\mF^{(i)}) = \kappa_i(\tX)$.
\end{proposition}

\begin{proof}
  From the preceding comments, we have that $\rank_+(\mF^{(i)}) \leq
  \kappa_i(\tX)$. Since $\rank_+(\tX) = \rank_+(\tG)$ from the canonical
  nnTD, we have that $\tG$ admits a nonnegative rank $R$ nnCPD
  \[
  \tG  =  \tD\times_1\mB^{(1)}\times_2\mB^{(2)}\times_3\mB^{(3)}.\\
  \]
  Substituting into the canonical nnTD, we two rank $R$ decompositions of $\tX$ given by
  \[
  \begin{array}{rcl}
    \tX&=&\tD\times_1\mA^{(1)}\times_2\mA^{(2)}\times_3\mA^{(3)}\\
    &=&\tD\times_1(\mF^{(1)}\mB^{(1)})\times_2(\mF^{(2)}\mB^{(2)})\times_3(\mF^{(3)}\mB^{(3)}).\\
  \end{array}
  \]
  Thus up to scaling and permutation, $\mF^{(i)}\mB^{(i)} = \mA^{(i)}$
  so that $\rank_+(\mF^{(i)}) \geq \rank_+(\mA^{(i)}) =
  \kappa_i(\tX)$.
\end{proof}

What Theorem \ref{thm: canonical preserve rank} and Proposition
\ref{prop: canonical} indicate is that a nonnegative canonical Tucker
decomposition gives just enough breathing room to provide existence of
a not rank degenerate Tucker as seen in
Figure~\ref{fig:nTDschematic}. In the diagram, we see the relation
that if a minimal nnTD also preserves the rank, then the minimal
multiranks are also the canonical multiranks. There are rank
preserving nnTDs that are neither minimal, nor canonical.  These can
be obtained, for instance, by taking a rank preserving canonical and
adding appropriate zeros.
As noted above, not every nonnegative canonical TD can preserve the
rank, and identifying the correct extensions of the minimal cones to
higher order cones to construct a rank preserving nnTD is a challenge.
In the next section, we compare the performance of different
nonnegative CANDELINC algorithms to demonstrate these challenges on
synthetic and real data.

\section{Numerical Experiments}
\label{sec:numerical}

Exact decompositions of tensors from real or experimental data are
typically not attainable, so we solve for approximations with
optimization problems. To find a nnCANDELINC approximation of a tensor
$X \in \mathbb{R}^{n_1, n_2, n_3}_+$, given a set of multirank
dimensions $[m_1,m_2,m_3]$ and a tensor rank $r$, a typical Frobenius
norm optimization problem is,
\begin{align}
  \label{eq:nnCANDELINC}
  \begin{array}{lll}
    &\underset{\mF^{(i)},\mB^{(i)}}{\operatorname{minimize}} & \normF{\tX-
\tI \times_1\mF^{(1)}\mB^{(1)}\times_2\mF^{(2)}\mB^{(2)}\times_3\mF^{(3)}\mB^{(3)} }^2 \\[0.5em]
    &\textrm{Subject to: }  & \mF^{(i)} \in \mathbb{R}^{n_i,m_i}_+ \text{ for } 1 \leq i \leq 3\\
    && \mB^{(i)} \in \mathbb{R}^{m_i, r}_+ \text{ for } 1 \leq i \leq 3\;.
  \end{array}
\end{align}

A similar optimization problem, without the nonnegativity constraints,
is suitable for CANDELINC.
To compute CANDELINC decompositions, theory informs us of two
procedures to find approximate decompositions using the readily
available tools of TD, CPD and SVD.
One procedure follows as: first compute a TD followed by a CPD on the
core.
Alternatively: first compute a CPD followed by SVDs on each of the CPD
factors.
Theoretically, both of these two procedures provide equally valid
CANDELINC decompositions under perfect conditions.
In practice the first procedure is often preferred as it is typically
less computationally expensive with the cheap dimension reduction via
Tucker compression before the more expensive CPD step.
In this section, we explore two procedures to compute nnCANDELINC
decompositions. We discuss the theory, benefits and drawbacks, and
demonstrate their performance on synthetic and real data.

\subsection{Algorithms and Scoring}
While a multitude of methods can be applied to the nnCANDELINC
optimization in Equation \ref{eq:nnCANDELINC}, we concentrate on
procedures that can be built using readily available tools. Namely, we
are interested in procedures that compute nnCANDElINC decompositions
using the sub-procedures: nnTD, nnCPD, and NMF.

Algorithm~\ref{alg:ntd-ncpd} computes a nnCANDELINC decomposition by
first nnTD compression, reducing the dimension, and then by nnCPD on
the resulting core. Much of the previously discussed theory informs us
of potential problems with this procedure in selecting various
multiranks and ranks for the nnTD and nnCPD dimensions.  Namely, it is
necessary, but not sufficient, for the nnTD to be computed with
dimensions greater or equal to the canonical multiranks, not the
minimal multiranks, to obtain a nnCANDELINC decomposition that
expresses the tensors rank.

\begin{algorithm}
  \caption{nnCANDELINC algorithm using nnTD and nnCPD.}
  \label{alg:ntd-ncpd}
  \begin{algorithmic}
    \Require $\mathcal{X} \in \mathbb{R}^{N_1, N_2, N_3}_+, m \in \mathbb{N}^3, r \in \mathbb{N}$
    \Ensure $\mF^{(i)} \in \mathbb{R}^{N_i, m_i}_+, \mB^{(i)} \in \mathbb{R}^{m_i,r}_+$
    \State $\tG, \mF^{(1)}, \mF^{(2)}, \mF^{(3)} \gets \operatorname{nnTD}(X, (m_1, m_2, m_3))$
    \State $\mB^{(1)}, \mB^{(2)}, \mB^{(3)} \gets \operatorname{nnCPD}(\tG, r)$
  \end{algorithmic}
\end{algorithm}

Alternatively, Algorithm~\ref{alg:ncpd-nmf} computes a nnCANDELINC
decomposition by first computing nnCPD, followed by NMF on each of the
nnCPD factors. Here, theory informs us that it is both necessary and
sufficient to use the canonical multiranks as the latent dimensions in
their respective NMFs.

\begin{algorithm}
  \caption{nnCANDELINC algorithm using nnCPD and NMF.}
  \label{alg:ncpd-nmf}
  \begin{algorithmic}
    \Require $\mathcal{X} \in \mathbb{R}^{N_1, N_2, N_3}_+, m \in \mathbb{N}^3, r \in \mathbb{N}$
    \Ensure $\mF^{(i)} \in \mathbb{R}^{N_i, m_i}_+, \mB^{(i)} \in \mathbb{R}^{m_i,r}_+$
    \State $\mA^{(1)}, \mA^{(2)}, \mA^{(3)} \gets \operatorname{nCPD}(X, r)$
    \State $\mF^{(1)}, \mB^{(1)} \gets \operatorname{NMF}(\mA^{(1)}, m_1)$
    \State $\mF^{(2)}, \mB^{(2)} \gets \operatorname{NMF}(\mA^{(2)}, m_2)$
    \State $\mF^{(3)}, \mB^{(3)} \gets \operatorname{NMF}(\mA^{(3)}, m_3)$
  \end{algorithmic}
\end{algorithm}

To demonstrate the performance of the Algorithms we utilize the
functions

\begin{center}
  \textbf{tensorly.decomposition.non\_negative\_tucker}
\end{center}

\noindent
and
\noindent

\begin{center}
\textbf{tensorly.decomposition.non\_negative\_parafac}
\end{center}

\noindent
from the freely available high-level API for tensor decomposition
methods in python, TensorLy~\cite{kossaifi2019tensorly}, and

\begin{center}
\textbf{sklearn.decomposition.NMF}
\end{center}

\noindent
from scikit-learn~\cite{scikit-learn}. In all experiments, random
initializations are used with each call and a constant 5000 iterations
are used for each sub-optimization to ensure reasonable convergence,
with no early termination criteria.

To evaluate the resulting decompositions of these algorithms we
utilize two different scores. The
congruence~\cite{tomasi2006comparison} between two rank one tensors,
$\mathcal{X}=a_1 \otimes b_1 \otimes c_1$ and $\mathcal{Y} = a_2
\otimes b_2 \otimes c_2$ is,
$$
\operatorname{cong}(\mathcal{X},\mathcal{Y}) = \cos(X,Y) =
\frac{a_1^\top \cdot a_2}{\Vert a_1 \Vert_2 \Vert a_2 \Vert_2} \cdot
\frac{b_1^\top \cdot b_2}{\Vert b_1 \Vert_2 \Vert b_2 \Vert_2} \cdot
\frac{c_1^\top \cdot c_2}{\Vert c_1 \Vert_2 \Vert c_2 \Vert_2}.
$$
The mean congruence of all rank one factors is relevant for two
rank $r$ tensors after the appropriate permutation of the factors is
applied to maximize the mean
congruence~\cite{battaglino2018practical}. We apply the mean
congruence to the appropriate products, $\mF^{(i)} \mB^{(i)}$, in the
nnCANDELINC decompositions. For a tensor $\mathcal{X}$, the Frobenius
norm is defined as the square root of the sum of the squares or
$||\mathcal{X}||_F = \sqrt{\sum_{i,j,k} \mathcal{X}_{i,j,k}^2 }$. The
relative reconstruction error (we call further relative error) of the
decomposition is the ratio of the Frobenius norm of the residual and
the Frobenius norm of the tensor or matrix. In addition to the mean
congruence, we utilize the relative reconstruction error to evaluate
the quality of the nnCANDELINC decompositions.

To evaluate these algorithms and relate them to theory we apply them
to both synthetic and real datasets.  Our first investigation of a
synthetic tensor highlights the importance of using the nonnegative
canonical multiranks and not the minimal multiranks.  We additionally
construct a large number nonnegative tensors with pre-determined
nonnegative canonical multiranks, and show the performance of each
algorithm at recovering the factors to evaluate them in more generic
situations.  For the first real dataset, we apply the nnCANDELINC
algorithms to extract the nCPD features of a well-known fluorescence
data that has been previously analyzed in the PhD Thesis of
Bro~\cite{Bro:1998:PhD}.
Next, we apply nnCANDELINC to a computer generated $3D$ dataset with
nonnegative rank deficient nCPD factors that represents a microphase
separation of block copolymers as a function of temperature and was
analyzed in~\cite{Alexandrov-etaL:2019}.

\subsection{Various Multiranks}

\begin{figure}
  \centering
  \begin{subfigure}[b]{0.45\linewidth}
    \centering
    \includegraphics[width=\textwidth]{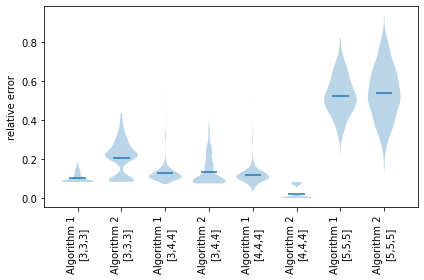}
    \caption{Violin plots of relative errors of decompositions.}
    \label{fig:minimal_canonical_error}
  \end{subfigure}
  \hfill
  \begin{subfigure}[b]{0.45\linewidth}
    \centering
    \includegraphics[width=\textwidth]{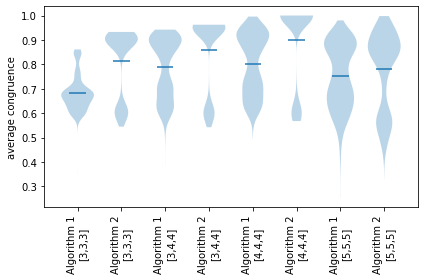}
    \caption{Violin plots of congruence scores of decompositions.}
    \label{fig:minimal_canonical_congruence}
  \end{subfigure}
  \caption{Violin plots and means of the results of nnCANDELINC
    decompositons on the tensor from Example~\ref{ex:equal} using the
    various nonnegative multiranks in Algorithms~\ref{alg:ntd-ncpd}
    and~\ref{alg:ncpd-nmf}.}
  \label{fig:minimal_canonical}
\end{figure}

We first investigate the efficacy of Algorithms~\ref{alg:ntd-ncpd} and
\ref{alg:ncpd-nmf} when various multiranks are used for the nnTD and
NMF dimensions. Example~\ref{ex:equal} provides an instance where the
nonnegative rank, minimal multiranks, and canonical multiranks are all
known. We evaluate nnCANDELINC on the proposed in the
Example~\ref{ex:equal} tensor with four different effective
nonnegative multiranks:
\begin{itemize}
\item $[3,3,3]$ which are less than the minimal multiranks of the
  tensor,
\item $[3,4,4]$ which are the minimal multiranks which is known not to
  preserve the rank to the core,
\item $[4,4,4]$ which are the canonical multiranks where it is
  feasible that the rank is preserved to the core,
\item $[5,5,5]$ which are greater than the latent dimensions needed
  everywhere.
\end{itemize}

To evaluate the algorithms we decompose the Example~\ref{ex:equal}
tensor 1000 times with each set of assumed multiranks, each starting
from random initial conditions.

Figure~\ref{fig:minimal_canonical} reports violin plots of the
relative errors and average congruence scores of the resulting
decompositions.
In Figure~\ref{fig:minimal_canonical_error} we see mild relative
errors for the three smallest assumed nonnegative multiranks $[3,3,3],
[3,4,4]$, and $[4,4,4]$.
The multiranks $[3,3,3]$ and $[3,4,4]$ are not sufficient to obtain a
rank revealing nnCANDELINC decomposition resulting in moderate
relative errors.
Using the canonical nonnegative multiranks, $[4,4,4]$ is expectedly
more successful with Algorithm~\ref{alg:ncpd-nmf}.
This is unsurprising since for both Algorithm~\ref{alg:ntd-ncpd}
and~\ref{alg:ncpd-nmf} the canonical multiranks are necessary to
obtain a low relative error, but for Algorithm~\ref{alg:ncpd-nmf} the
use of the canonical multiranks is both necessary and sufficient.
Using an assumed multirank of $[5,5,5]$ leads to disastrous
performance with both algorithms, on the surface this is surprising
since the problem has more degrees of freedom than necessary to
perfectly reconstruct the tensor.
Figure~\ref{fig:minimal_canonical_congruence} reflects the mean and
standard deviations of the average congruences for each scheme.
There is little correlation between the relative error performances
and the congruence scores.
By the congruence score measure, Algorithm~\ref{alg:ncpd-nmf} shows
better scores than Algorithm~\ref{alg:ntd-ncpd} on this tensor.

\subsection{Randomly Generated Tensors}

\begin{figure}
  \centering
  \begin{subfigure}[b]{0.45\linewidth}
    \centering
    \includegraphics[width=\textwidth]{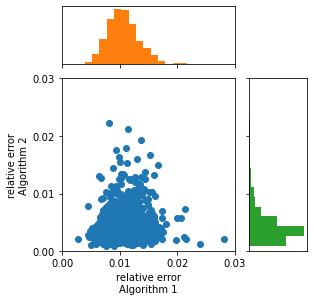}
    \caption{Relative errors of decompositions of synthetic dataset.}
    \label{fig:minimal_canonical_error:b}
  \end{subfigure}
  \hfill
  \begin{subfigure}[b]{0.45\linewidth}
    \centering
    \includegraphics[width=\textwidth]{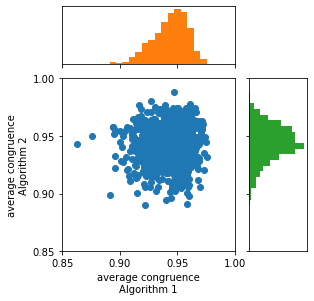}
    \caption{Average congruence scores of decompositions of synthetic dataset.}
    \label{fig:minimal_canonical_congruence:b}
  \end{subfigure}
  \caption{Scatterplots and projected histograms from results of
    nnCANDELINC decompositons using Algorithms~\ref{alg:ntd-ncpd}
    and~\ref{alg:ncpd-nmf} on randomly generated synthetic
    dataset. Each point corresponds to a tensor whose location
    indicates how well each algorithm performed.}
    \label{fig:canonical}
\end{figure}

To evaluate the performance of the nnCANDELINC Algorithms on a more
varied dataset, we randomly generate synthetic tensors that have a
unique nnCPD with nonnegative rank deficient factors and where the
canonical multiranks are known. To do this, we randomly generate our
nnCPD factors as a product of two smaller nonnegative matrices, and
confirm that they satisfy the suppositions of Kruskal's theorem
\cite{kruskal1977three}. For specified dimensions $N_1, N_2, N_3$,
nonnegative canonical multiranks $r_1, r_2, r_3$, and rank $r$, we
randomly sample from a uniform distribution the factors $F^{(i)} \in
\REAL_+^{N_i \times r_i}$ and $B^{(i)} \in \REAL_+^{r_i \times r}$ to
construct the decomposition $\tX = \tI \times_1 (F^{(1)} B^{(1)})
\times_2 (F^{(1)} B^{(2)}) \times_3 (F^{(1)} B^{(3)})$. We ensure that
factors satisfy the Kruskal rank criteria,
$\mbox{Kruskal-rank}(F^{(1)} B^{(1)}) + \mbox{Kruskal-rank}(F^{(2)}
B^{(2)}) + \mbox{Kruskal-rank}(F^{(3)} B^{(3)}) \leq 2r+2$.

Here we report the results of 1000 randomly generated tensors, each of
size $N_1 = N_2 = N_3 = 40$ with a nonnegative rank of $5$ and
canonical multiranks of $[3,4,5]$.
Figure~\ref{fig:canonical} depicts scatterplots of the resulting
relative errors and congruence scores using each algorithm for each of
the 1000 tensors.  Both algorithms yield a low relative error and high
average congruence score with Algorithm~\ref{alg:ncpd-nmf}
demonstrating better scores than Algorithm~\ref{alg:ntd-ncpd}.

\subsection{Fluorescence Data Decompositions}
\begin{figure}[tbp]
  \centering
  \begin{subfigure}[b]{0.45\linewidth}
    \centering
    \includegraphics[width=\textwidth]{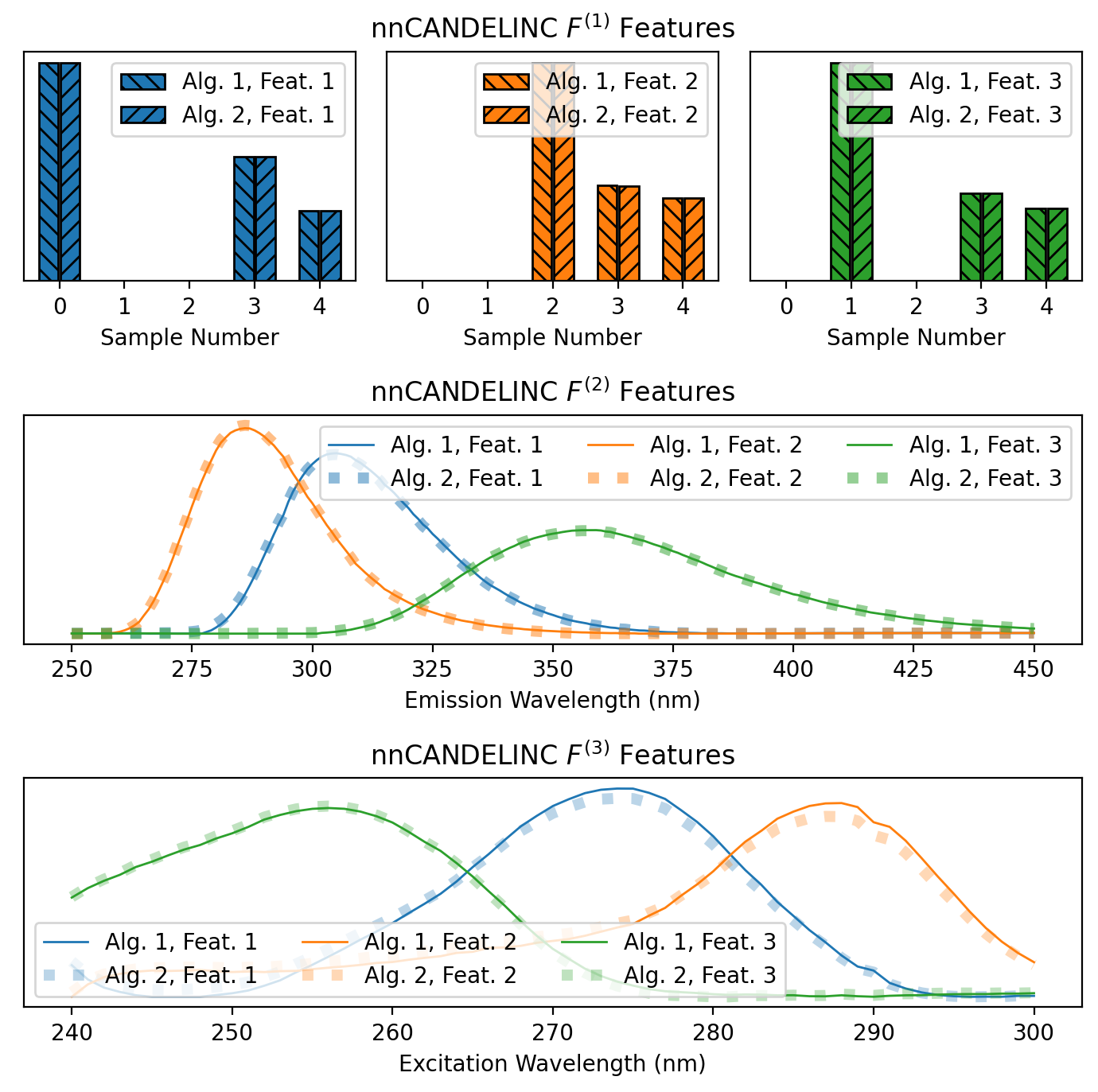}
    \caption{Comparison of features $F^{(i)}$ as extracted from the
      fluorescence tensor using both Algorithms~\ref{alg:ntd-ncpd}
      and~\ref{alg:ncpd-nmf}}
    \label{fig:fluorF}
  \end{subfigure}
  \hfill
  \begin{subfigure}[b]{0.45\linewidth}
    \centering
    \includegraphics[width=\textwidth]{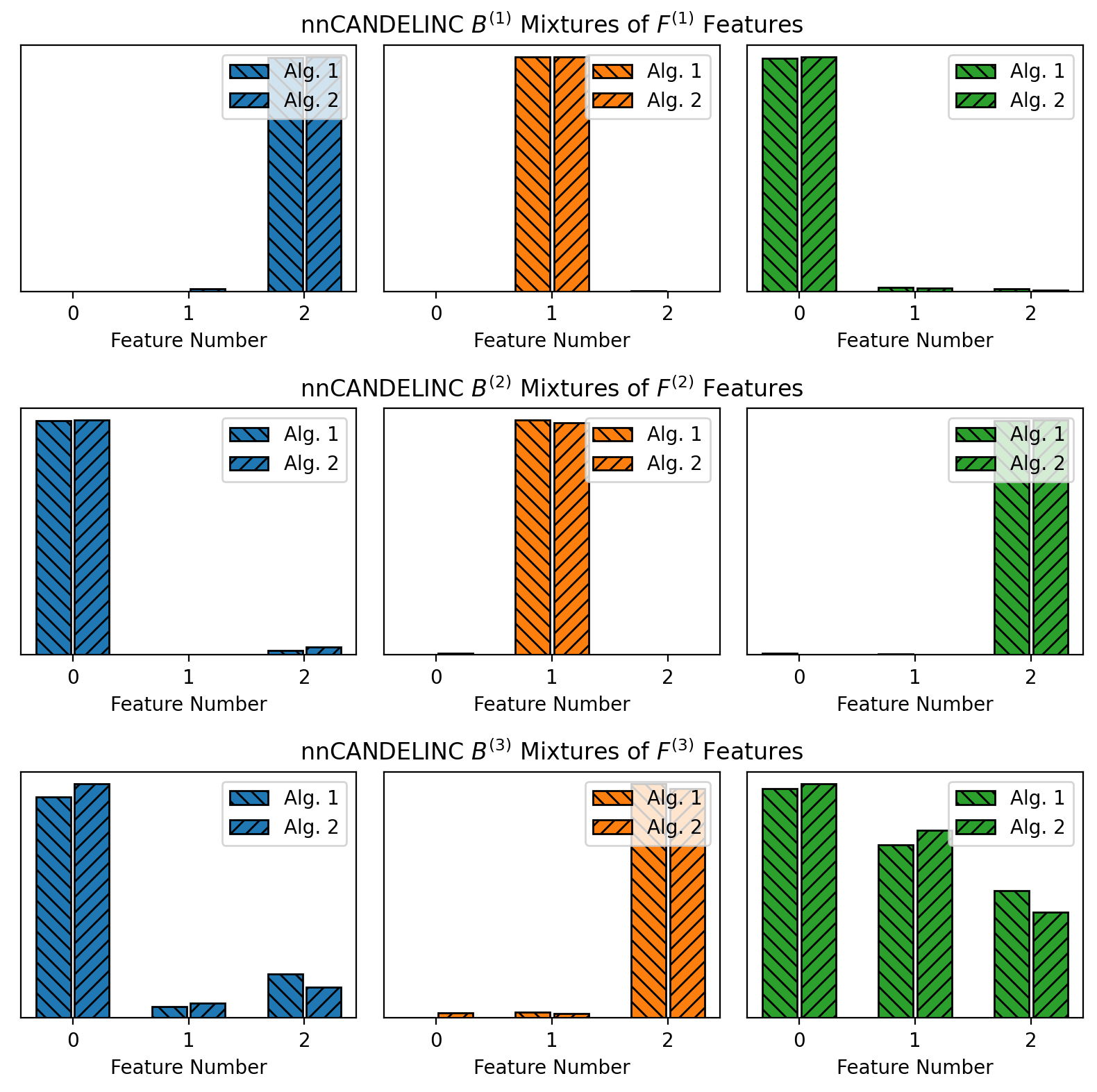}
    \caption{Comparison of activations $B^{(i)}$ of the features
      $F^{(i)}$ from the fluorescence tensor using both
      Algorithms~\ref{alg:ntd-ncpd} and~\ref{alg:ncpd-nmf}}
    \label{fig:FluorB}
  \end{subfigure}
  \caption{ Comparisons of the features, $\mF^{(i)}$, and mixings,
    $\mB^{(i)}$, obtained from nnCANDELINC
    Algorithms~\ref{alg:ntd-ncpd} and~\ref{alg:ncpd-nmf} on the
    experimental fluorescence data of size ($5 \times 61 \times 201$)
    with nonnegative multiranks $[3,3,3]$, and nonnegative rank
    $r$=3. }
  \label{Flour:nCPDk}
\end{figure}

The experimental fluorescence dataset includes five samples, each with
different amounts of amino acids of three types: tyrosine, tryptophan
and phenylalanine dissolved in buffered water.  We consider this data
to demonstrate the different algorithm performances when there is no
strong linear dependence of the factors.
The fluorescence in these samples has been excited by UV irradiation
at wavelengths, $\lambda\in (240-300\textrm{nm})$.
The UV-emission was measured by the spectrofluorometer at wavelengths
$\lambda\in[250,450]\textrm{nm}$ by sampling at 1 nm intervals.
The experimental data formed a $3D$ array with size $5\times 61\times
201$.
If we assume that each amino acid gives a nonnegative linear
contribution to the fluorescence data-tensor, than the measured
fluorescence, i.e., the output, $\tX$, is three-linear, and its
components, $\tX_{i, j, k}$ are,
\begin{align*}
  \tX_{i,j,k} = \sum^{r}_{n=1}\mA^{(1)}_{i,n}\,\mA^{(2)}_{j,n}\,\mA^{(3)}_{k,n}+\epsilon_{i,j,k}.
\end{align*}
Here, 
$\mA^{(1)}_{i,n}$ is linearly related to the concentration of the
$n^{th}$ fluorophore dissolved in the $i^{th}$ sample;
$A^{(2)}_{j,n}$ to the relative emission of $n^{th}$ fluorophore at
wavelength $\lambda_j$; 
$A^{(3)}_{k,n}$ to the relative amount of UV light absorbed by
$n^{th}$ fluorophore at excitation $\lambda_k$, and $\epsilon_{i,j,k}$
denotes the error.
Although the above formula represents an ideal physical situation, it
has been shown that for small concentrations of amino acids it is a
valid approximation ~\cite{Bro:1998:PhD}.
Here, we apply the nnCANDELINC Algorithms described in the previous
sections, compare their results, and validate that the final
decompositions coincides with the previously well-known results.

The scaled and appropriately permuted results presented in
Figure~\ref{Flour:nCPDk} show minor differences between the resulting
decompositions obtained from Algorithm~\ref{alg:ntd-ncpd} with 2.79\%
relative error, and Algorithm~\ref{alg:ncpd-nmf} with
2.51\%. Figure~\ref{fig:fluorF} depicts the features extracted along
the sample, emission, and excitation axes. We see virtually no
difference in the sample and emission extracted features, with only
slight deviations occurring in the excitation features. Similarly in
Figure~\ref{fig:FluorB} the mixtures of the sample and emission
features are virtually identical, while there are slight deviations in
the mixtures of the excitation features between the two algorithms.  A
comparison with previously extracted features from the same
data~\cite{Bro:1998:PhD} confirms that both nnCANDELINC algorithms are
producing correct results.
It is also worth mentioning that the utilization of the nonnegative TD
in Algorithm~\ref{alg:ntd-ncpd} does not results in a superdiagonal
core-tensor $\tG$, and the products of the final factors of both
nnCANDELINC algorithms are indistinguishable from those obtained by a
direct application of CPD~\cite{Bro:1998:PhD}.  

\subsection{Decomposition of data generated by physics-based computer simulations}

\begin{figure}[tbp]
  \centering
  \begin{subfigure}[b]{0.45\linewidth}
    \centering
    \includegraphics[width=\textwidth]{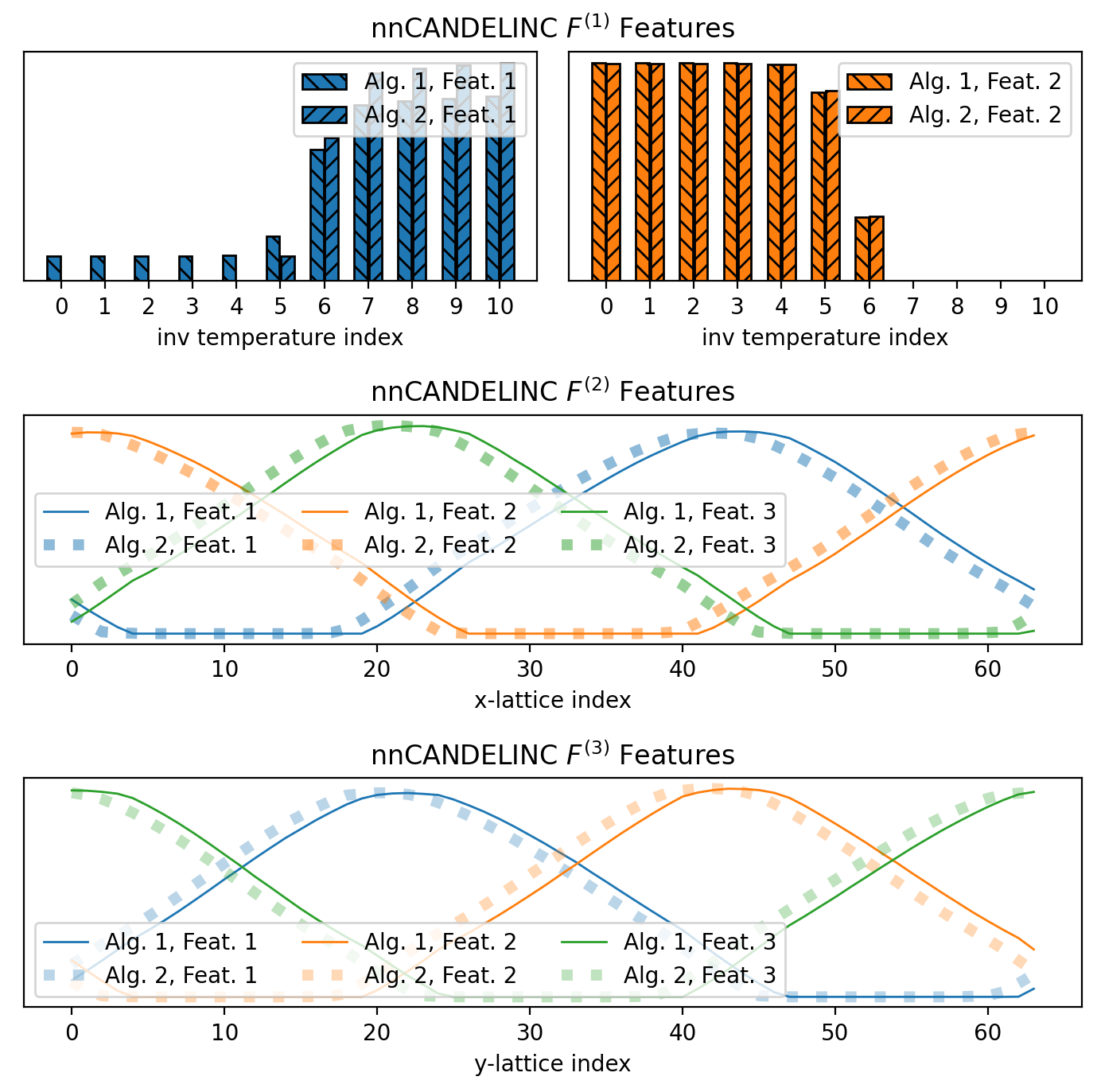}
    \caption{Comparison of features $F^{(i)}$ as extracted from the
      copoloymers tensor using both Algorithms~\ref{alg:ntd-ncpd}
      and~\ref{alg:ncpd-nmf}}
    \label{fig:polyF}
  \end{subfigure}
  \hfill
  \begin{subfigure}[b]{0.45\linewidth}
    \centering
    \includegraphics[width=\textwidth]{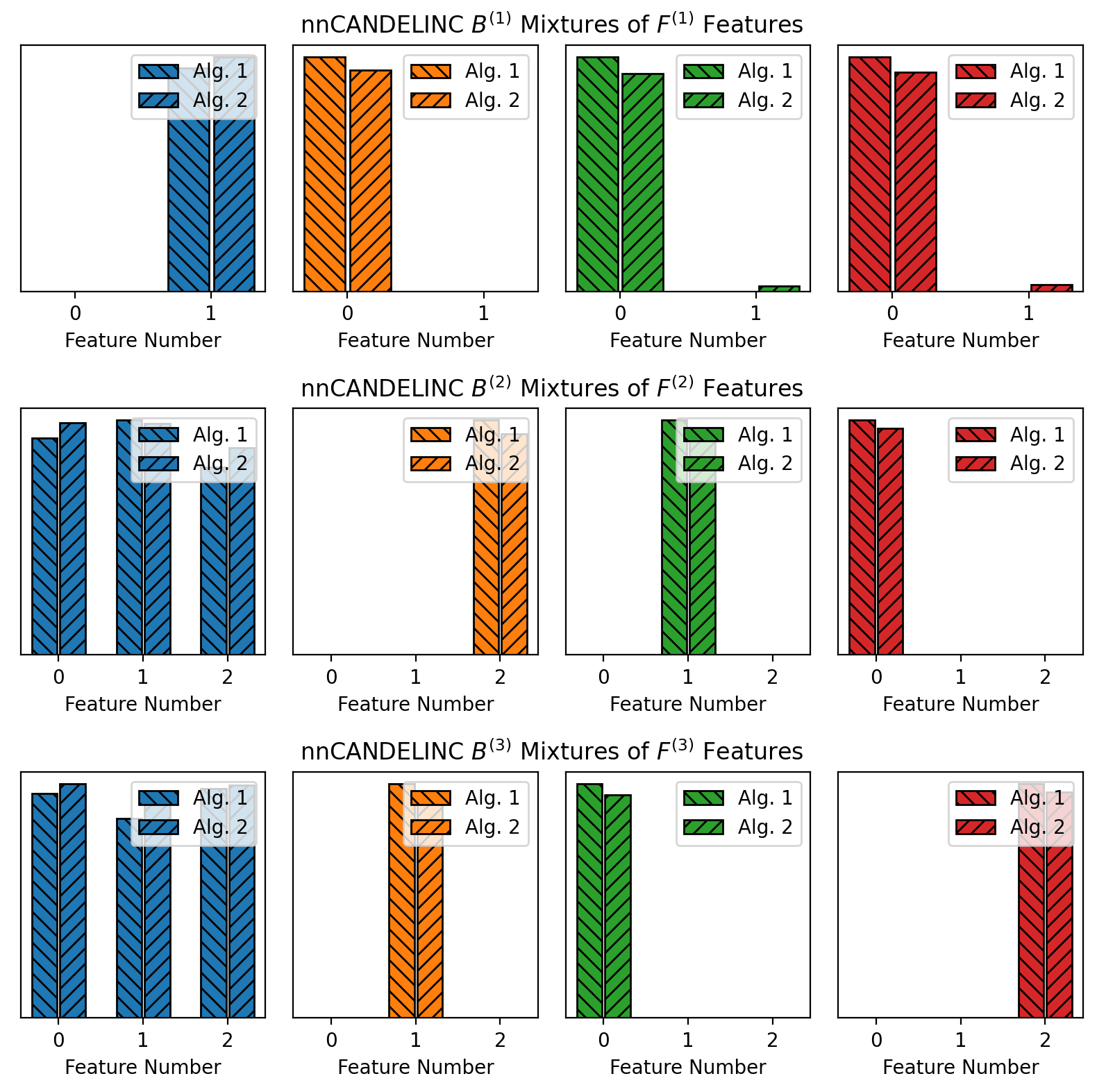}
    \caption{Comparison of activations $B^{(i)}$ of the features
      $F^{(i)}$ from the copolymers tensor using both
      Algorithms~\ref{alg:ntd-ncpd} and~\ref{alg:ncpd-nmf}}
    \label{fig:polyB}
  \end{subfigure}
  \caption{
    Performances of nnCANDELINC Algorithms~\ref{alg:ntd-ncpd}
    and~\ref{alg:ncpd-nmf} on computer generated data with size ($11
    \times 64 \times 64$), nonnegative multiranks $[2,3,3]$ and rank
    $r=4$, representing phase separation with temperature in a system
    of copolymers.  }
  \label{Phase:nCPDk}
\end{figure}

Here, we use nnCANDELINC to analyze a $3D$ data-tensor describing
phase separation in a system of blocks copolymers whose evolution with
temperature has been introduced and analyzed in a previous
work~\cite{Alexandrov-etaL:2019}.
We chose this system because of the natural nonnegativity of the data,
the already known nonnegative rank, $r$ = 4, and the fact that the
extracted factors have a rank deficiency demonstrated in the previous
analysis.

The multivariate function describing the phase separation is the order
parameter of the system, $\Delta(T, f_A, x, y)$, which in this case is
a function of: (a) temperature, $T$, (b) length $f_A$ of the A-type
blocks, and (c) the spatial coordinates, $(x, y)$, of the
$2$-dimensional $64 \times 64$ lattice-space of the system.
The order parameter, $\Delta(T, f_A, x, y)$, is simply the spatial
density of the A-type blocks on the lattice, and therefore the data is
inherently nonnegative. For A-type blocks with a fixed length, $f_A$,
the order parameter is represented by $3$-dimensional data: $\Delta(T,
f_A, x, y)\equiv \Delta(T, x, y)$, and the tensor $\Delta_{n, m, l}$
that we analyze here has size $11\times 64 \times 64 $.

The nonnegative ranks $r_{i}$ of each unfolding of the tensor
$\Delta(T, x, y)$ has been previously estimated
\cite{Alexandrov-etaL:2019} and the nonnegative multirank has been
determined to be,
$\multirank_+(\tX) = [2,3,3]$. 
With this nonnegative minimal multirank we applied both nnCANDELINC
algorithms to $\tX$ and compare the results.

In Figure \ref{Phase:nCPDk} we present the components of factors
$\mF^{(i)}$ and $\mB^{(i)}$ from both Algorithm~\ref{alg:ntd-ncpd}
with a relative error of 11.29\%, and Algorithm~\ref{alg:ncpd-nmf}
with an error of 10.02\%.
In Figure \ref{fig:polyF} the extracted features from the algorithms
vary slightly, along the temperature axis we see relative shifts
between the feature extracted by the algorithms, with similar shifts
seen in the x-lattice and y-lattice axes. These shifts result in
slight differences of the mixing of these features seen in
Figure~\ref{fig:polyB}. The nonnegative rank deficiencies become clear
with the found combinations of features to represent the four rank one
tensors needed for an nnCPD.  

\appendix
\section{Notation and Operations}

Here we list the precise operations used throughout the paper. A
useful operation often used is the multiplication of a tensor by a
matrix along a specific dimension, or \emph{$n$-mode multiplication}.
\begin{definition}
  \label{def:n-mode}
  The \emph{$1$-mode multiplication} between a tensor
  $\tX\in\REAL^{N_1\times N_2\times N_3}$ and a matrix
  $\mA\in\REAL^{M\times N_1}$ is defined as
  \begin{align*}
    (\tX\times_1\mA)_{i,j,k} = \sum_{l=1}^{N_1}\tX_{l,j,k}\mA_{i,l}.
  \end{align*}
  We define the $2$-mode and $3$-mode multiplication analogously.
\end{definition}
For $i\neq j$ mode multiplications are commutative:
$(\tX\times_i\mA^{(i)})\times_j\mA^{(j)}=(\tX\times_j\mA^{(j)})\times_i\mA^{(i)}$,
and a matrix multiplication can be distributed through mode
multiplication: $\tX\times_i\mA\mB=(\tX\times_i\mB)\times_i\mA$.
\begin{definition}
  \label{def:fiber}
  A \emph{mode-$i$ tensor fiber} of $\tX$ is a one dimensional vector
  obtained by fixing all but the $i^{th}$ index in the tensor.
  We let $\tX_{:,n,m}$, $\tX_{n,:,m}$, $\tX_{n,m,:}$ denote the
  $n^{th}, m^{th}$ mode-$1$, mode-$2$ and mode-$3$ tensor fibers,
  respectively. For $i=1,2,3$, $\matricization_i(\tX)$ denotes an \emph{$i$-mode
  unfolding}, which rearranges all the mode-$i$ fibers of a tensor
  into columns of a $N_i\times N_jN_k$ matrix, for $i\neq j\neq k$.
\end{definition}
\begin{definition} 
  \label{def:unfolding}
  For $i=1,2,3$, $\matricization_i(\tX)$ denotes an \emph{$i$-mode
    unfolding}, which rearranges all the mode-$i$ fibers of a tensor
  into columns of a $N_i\times N_jN_k$ matrix, for $i\neq j\neq k$.
\end{definition}
Each unfolding has an inverse mapping, which rearranges the columns of
a matrix as fibers of a tensor.
Consider the tensor
$\tX=\tY\times_1\mA^{(1)}\times_2\mA^{(2)}\times_3\mA^{(3)}$.
A particularly useful relation between unfoldings and mode
multiplications is
\begin{align}
  \matricization_i(\tX) = A^{(i)} \matricization_i(\tY \times_j
  A^{(j)} \times_k A^{(k)}),
  \label{eq:unfolding:property}
\end{align}
for $i\neq j\neq k$.

\section{Basics of NMF}
Nonnegative matrix factorization (NMF) decomposes a nonnegative matrix
$\mV\in\REAL^{N\times M}_+$, into a product of two nonnegative
matrices $\mW\in\REAL^{N\times r}_+$ and $\mH\in\REAL^{r\times M}$.
The geometric interpretation of a nonnegative decomposition, $V=WH$,
is that that each column of $V$ is a conic combination of the columns
of $W$.
With this geometric interpretation, computing an NMF is identical to searching for a
polyhedral cone $C$ which contains the columns of $\mV$, and is
contained in the nonnegative orthant, $V \subset C \subset \REAL^N_+$.
Of particular interest are cones with a minimum number of extreme
rays, which correspond to the nonnegative rank.

\begin{definition}
  The nonnegative rank of a matrix is defined as
  \begin{align*}
    \rank_+{(\mV)} := \min\left\{
    r~\bigg\vert\,
    \mV=\sum_{n=1}^r\vw_{n}\ccirc\vh_{n},\,
    \vw_n\geq 0,\,
    \vh_n\geq 0 
    \right\}.
  \end{align*}
\end{definition}
If $\rank_+(\mV) = r$ then there is a set of $r$ nonnegative extreme
rays $\{w_1, \hdots, w_r \}$ such that every column of $\mV$ is a
conic combination of these extreme rays.
When $\{w_1, \hdots, w_r \}$ are assembled into the nonnegative matrix
$\mW$, and the conic combinations are specified by a nonnegative
matrix $H$, this corresponds to the nonnegative matrix factorization
$\mV = \mW \mH$.

The nonnegative rank of a matrix has several well-known properties.
For example, if $\mV$ is an $(N_1\times N_2)$-sized matrix, then
$\rank(\mV)\leq\rank_+(\mV)\leq\min(N_1,N_2)$~\cite{Cohen-Rothblum-RUTCOR:1993}.
A case illustrating the inequality between rank and nonnegative rank
can be seen in the following Example \ref{ex:mrank_mrank+}, which is
mentioned in~\cite{Cohen-Rothblum-RUTCOR:1993} as a private
communication from H.~Robbins.

\begin{example}\label{ex:mrank_mrank+}
  Consider the nonnegative matrix:
  \begin{align}
    \mV = 
    \begin{bmatrix} 
      1 &\  1 &\  0 &\  0\\ 
      1 &\  0 &\  1 &\  0\\ 
      0 &\  1 &\  0 &\  1\\ 
      0 &\  0 &\  1 &\  1\\ 
    \end{bmatrix}
    \label{eq:matV:nnCPD}
  \end{align}
  and note that $v_1 + v_4 = v_2 + v_3$ where $v_i$ is the $i^{th}$
  column of V.
  This linear dependence between the columns proves that $\rank(V) =
  3$. Also it was proved in ~\cite{Cohen-Rothblum-RUTCOR:1993} that
  the $\rank_+(V) = 4$.
  This example demonstrates a case when $\rank(V) < \rank_+(V)$.
  \ENDPROOF
\end{example}

In general, computing the nonnegative rank of a nonnegative matrix
$\mV\in\REAL^{N_1 \times N_2}$ is an NP-hard
problem~\cite{Cichocki:Zdunek:Phan:Amari:2009}, and even providing a
reliable estimate can be quite hard.


\section*{Acknowledgments}
This work was supported by the LDRD program of Los Alamos National
Laboratory under project number 20190020DR.
Los Alamos National Laboratory is operated by Triad National Security,
LLC, for the National Nuclear Security Administration of
U.S. Department of Energy (Contract No. 89233218CNA000001).
We would also like to thank the anonymous reviewer for their helpful
comments and suggestions.
This study does not have any conflicts to disclose.



\begin{thebibliography}{10}

\bibitem{Big-Data}
Hickey A. 2019. \textit{Zettabytes of data hog up space and resources}.
\newblock

\bibitem{franke2016statistical}
Franke B., Plante J.F., Roscher R., Lee E.A., Smyth C., Hatefi A., et~al.
\newblock Statistical inference, learning and models in big data.
\newblock International Statistical Review. 2016;{\bf 84}(3):371--389.

\bibitem{Kolda-Bader:2009}
Kolda T., and Bader B.
\newblock Tensor decompositions and applications.
\newblock SIAM Review. 2009;{\bf 51}(3):455--500.

\bibitem{cichocki2017tensor}
Cichocki A., Phan A.H., Zhao Q., Lee N., Oseledets I., Sugiyama M., et~al.
\newblock Tensor networks for dimensionality reduction and large-scale
  optimization: {P}art~2 applications and future perspectives.
\newblock Foundations and Trends{\textregistered} in Machine Learning.
  2017;{\bf 9}(6):431--673.

\bibitem{oseledets2011tensor}
Oseledets I.V.
\newblock Tensor-train decomposition.
\newblock SIAM Journal on Scientific Computing. 2011;{\bf 33}(5):2295--2317.

\bibitem{vervliet2018compressed}
Vervliet N. 2018. \textit{Compressed sensing approaches to large-scale tensor
  decompositions}.
\newblock . KU Leuven.
\newblock (PhD Thesis).

\bibitem{Tucker:1966}
Tucker L.R.
\newblock Some mathematical notes on three-mode factor analysis.
\newblock Psychometrika. 1966;{\bf 31}(3):279--311.

\bibitem{Hackbusch:2012}
Hackbusch W.
\newblock Tensor spaces and numerical tensor calculus. vol.~42 of Springer
  series in computational mathematics.
\newblock Heidelberg: Springer; 2012.

\bibitem{Hitchcock:1927}
Hitchcock F.L.
\newblock The expression of a tensor or a polyadic as a sum of products.
\newblock Journal of Mathematics and Physics. 1927;{\bf 6}:164--189.

\bibitem{Harshman:1970}
Harshman R.A.
\newblock Foundation of the {PARAFAC} procedure: {M}odels and conditions for an
  ``explanatory'' multi-modal factor analysis.
\newblock In: UCLA Working Papers in Phonetics. vol.~16. University Microfilms,
  Ann Arbor, Michigan, No.~10,085; 1970. p. 1--84.

\bibitem{haastad1989tensor}
H{\aa}stad J.
\newblock Tensor rank is NP-Complete.
\newblock In: International Colloquium on Automata, Languages, and Programming.
  Springer; 1989. p. 451--460.

\bibitem{de2008tensor}
De~Silva V., and Lim L.H.
\newblock Tensor rank and the ill-posedness of the best low-rank approximation
  problem.
\newblock SIAM Journal on Matrix Analysis and Applications. 2008;{\bf
  30}(3):1084--1127.

\bibitem{everett2013introduction}
Everett B.
\newblock An introduction to latent variable models.
\newblock Springer Science \& Business Media; 2013.

\bibitem{lee1999learning}
Lee D.D, and Seung H.S.
\newblock Learning the parts of objects by non-negative matrix factorization.
\newblock Nature. 1999;{\bf 401}(6755):788.

\bibitem{Cichocki:Zdunek:Phan:Amari:2009}
Cichocki A., Zdunek R., Phan A.H., and Amari Si.
\newblock Nonnegative Matrix and Tensor Factorizations: Applications to
  Exploratory Multi-way Data Analysis and Blind Source Separation.
\newblock Wiley Publishing; 2009.

\bibitem{carroll1980candelinc}
Carroll J.D., Pruzansky S., and Kruskal J.B.
\newblock CANDELINC: A general approach to multidimensional analysis of
  many-way arrays with linear constraints on parameters.
\newblock Psychometrika. 1980;{\bf 45}(1):3--24.

\bibitem{bro1998improving}
Bro R., and Andersson C.A.
\newblock Improving the speed of multiway algorithms: Part II: Compression.
\newblock Chemometrics and intelligent laboratory systems. 1998;{\bf
  42}(1-2):105--113.

\bibitem{Bro:Harshman:Sidiropoulos:Lundy:2009}
Bro R., Harshman R.A., Sidiropoulos N.D., and Lundy M.E.
\newblock Modeling multi-way data with linearly dependent loadings.
\newblock Journal of Chemometrics: A Journal of the Chemometrics Society.
  2009;{\bf 23}(7-8):324--340.

\bibitem{wei2019compressed}
Wei W., Le X., Xiaofei Z., and Jianfeng L.
\newblock Compressed Sensing PARALIND Decomposition-Based Coherent Angle
  Estimation for 
\newblock Wireless Communications and Mobile Computing. 2019;{\bf 2019}.

\bibitem{Cohen:2016:PhD}
Cohen J. 2016. \textit{Environmental multiway data mining}.
\newblock . Universit\'e de Grenobles Alpes, France.
\newblock (Ph.D.~Thesis).

\bibitem{qi2018very}
Qi Y.
\newblock A Very Brief Introduction to Nonnegative Tensors from the Geometric
  Viewpoint.
\newblock Mathematics. 2018;{\bf 6}(11):230.

\bibitem{Qi-Comon-Lim:2016}
Qi Y., Comon P., and Lim L.H.
\newblock Semialgebraic geometry of nonegative tensor rank.
\newblock SIAM Journal of Matrix Analysis and Applications. 2016;{\bf
  37}(4):1556--1580.

\bibitem{landsberg2012tensors}
Landsberg J.M.
\newblock Tensors: geometry and applications.
\newblock Representation theory. 2012;{\bf 381}(402):3.

\bibitem{comon2014tensors}
Comon P.
\newblock Tensors: a brief introduction.
\newblock IEEE Signal Processing Magazine. 2014;{\bf 31}(3):44--53.

\bibitem{Bertsimas-Tsitsiklis:1997}
Bertsimas D., and Tsitsiklis J.N.
\newblock Introduction to linear optimization. vol.~6.
\newblock Athena Scientific Belmont, MA; 1997.

\bibitem{Gillis:2012}
Gillis N.
\newblock Sparse and unique nonnegative matrix factorization through data
  preprocessing.
\newblock Journal of Machine Learning Research. 2012;{\bf 13}(Nov):3349--3386.

\bibitem{Kruskal:1977}
Kruskal J.B.
\newblock Three-ways arrays: rank and uniqueness of trilinear decompositions,
  with application to arithmetic complexity and statistics.
\newblock Linear Algebra and Applications. 1977;{\bf 18}(2):95--138.

\bibitem{Domanov-DeLathauwer:2013}
Domanov I., and De~Lathauwer L.
\newblock On the uniqueness of the canonical polyadic decomposition of
  third-order tensors --- {P}art {I}: {B}asic results and uniqueness of one
  factor matrix.
\newblock SIAM Journal on Matrix Analysis and Applications. 2013;{\bf
  34}(3):855--875.

\bibitem{Cohen-Comon-Gillis:2017}
Cohen J.E., Comon P., and Gillis N.
\newblock Some theory on non-negative {T}ucker decomposition.
\newblock In: International Conference on Latent Variable Analysis and Signal
  Separation. Springer; 2017. p. 152--161.

\bibitem{kossaifi2019tensorly}
Kossaifi J., Panagakis Y., Anandkumar A., and Pantic M.
\newblock Tensorly: Tensor learning in python.
\newblock The Journal of Machine Learning Research. 2019;{\bf 20}(1):925--930.

\bibitem{scikit-learn}
Pedregosa F., Varoquaux G., Gramfort A., Michel V., Thirion B., Grisel O., et~al.
\newblock Scikit-learn: Machine Learning in {P}ython.
\newblock Journal of Machine Learning Research. 2011;{\bf 12}:2825--2830.

\bibitem{tomasi2006comparison}
Tomasi G., and Bro R.
\newblock A comparison of algorithms for fitting the PARAFAC model.
\newblock Computational Statistics \& Data Analysis. 2006;{\bf
  50}(7):1700--1734.

\bibitem{battaglino2018practical}
Battaglino C., Ballard G., and Kolda T.G.
\newblock A practical randomized CP tensor decomposition.
\newblock SIAM Journal on Matrix Analysis and Applications. 2018;{\bf
  39}(2):876--901.

\bibitem{Bro:1998:PhD}
Bro R. 1998. \textit{Multi-way analysis in the food industry. {M}odels,
  algorithms, and applications}.
\newblock . Royal Veterinary and Agricultural University, Denmark.
\newblock (Ph.D.~Thesis).

\bibitem{Alexandrov-etaL:2019}
Alexandrov B.S., Stanev V.G., Vesselinov V.V., and Rasmussen K.{\O}.
\newblock Nonnegative tensor decomposition with custom clustering for
  microphase separation of block copolymers.
\newblock Statistical Analysis and Data Mining: The ASA Data Science Journal.
  2019;.

\bibitem{kruskal1977three}
Kruskal J.B.
\newblock Three-way arrays: rank and uniqueness of trilinear decompositions,
  with application to arithmetic complexity and statistics.
\newblock Linear algebra and its applications. 1977;{\bf 18}(2):95--138.

\bibitem{Cohen-Rothblum-RUTCOR:1993}
Cohen J.E., Rothblum U.G., and {RUTCOR--Rutgers~Center~for~Operational~Research}.
\newblock Nonnegative ranks, decompositions and factorizations of nonnegative
  matrices.
\newblock Linear Algebra and its Applications. 1993;{\bf 190}:149--168.

\end{thebibliography}

\end{document}